\documentclass[12pt]{article}
\usepackage{amsmath,amssymb}
\providecommand{\MR}{\relax\ifhmode\unskip\space\fi MR }

\providecommand{\href}[2]{#2}

\newcommand{\qed}{{\unskip\nobreak\hfil\penalty50\hskip2em\vadjust{}
       \nobreak\hfil$\Box$\parfillskip=0pt\finalhyphendemerits=0\par}}

\newtheorem{theorem}{Theorem}[section]
\newtheorem{lemma}{Lemma}[section]

\newtheorem{cor}{Corollary}
\newtheorem{prop}{Proposition}

\renewcommand{\>}{\rangle}

\renewcommand{\a}{\alpha}
\renewcommand{\b}{\beta}
\newcommand{\bs}{\bigskip}

\newcommand{\ds}{\displaystyle}
\renewcommand{\d}{\delta}
\newcommand{\D}{\Delta}
\newcommand{\e}{\varepsilon}
\newcommand{\eps}{\epsilon}
\newcommand{\g}{\gamma}
\newcommand{\G}{\Gamma}
\renewcommand{\i}{\infty}

\renewcommand{\o}{\omega}

\newcommand{{\z}}{\mathbb Z}
\newcommand{\R}{\mathbb R}
\newcommand{\N}{\mathbb N}

\newcommand{\ue}{u^{\epsilon}}

\newcommand{\macp}{\beta}
\newcommand{\micp}{\alpha}

\renewcommand{\a}{\alpha}
\renewcommand{\b}{\beta}
\renewcommand{\d}{\delta}

\renewcommand{\o}{\omega}

\renewcommand{\i}{\infty}

\newcommand{\cH}{{\mathcal H}}
\newcommand{\cF}{{\mathcal F}}
\newcommand{\cG}{{\mathcal G}}

\newcommand{\parsize}{\big\vert \log \epsilon \big\vert^{-1}}
\newcommand{\sparsize}{\big\vert \log \epsilon \big\vert^{-2}}

\begin{document}
\setlength{\baselineskip}{16pt}

\begin{center}
{\Large\bf The kinetic limit of a system of coagulating planar
Brownian particles }\bs\bs\\
{\sc Alan Hammond$^1$}  \\   Department of
Statistics\\ University of California\\
Berkeley, California 93720 \bs \\
{\sc Fraydoun Rezakhanlou\footnote{Research
supported in part by NSF grant DMS0307021}}
  \\ Department of
Mathematics\\ University of California\\
Berkeley, California 93720--3830
\end{center}

\begin{abstract}We study a model of mass-bearing coagulating planar Brownian particles. 
Coagulation occurs when two particles are within a distance of order
$\e$. We assume that the initial number of particles $N$ is of order 
$|\log \e|$. Under suitable assumptions of the initial distribution of particles
and the microscopic coagulation propensities, we show that 
the macroscopic particle densities satisfy a Smoluchowski-type equation. 
\end{abstract}

\begin{section}{Introduction}

A colloid consists of a large number of small particles that are
suspended in an environment of far smaller and more numerous
molecules. Large numbers of molecules bombard each particle, and
random fluctuations among these collisions tend to give rise to an
Ornstein-Uhlenbeck 
motion of the particle, in which its velocity is forced by a
Brownian motion, with a drag force acting in the direction opposite to
its velocity. On a long time scale, the colloidal particles move
according to Brownian motions, because an Ornstein-Uhlenbeck process
approximates such a motion over a long period of time.
The particles of a colloid may also be liable to interact. In \cite{HR}, we
studied a model of a colloid in which this means of interaction took
the form of a coagulation, this reaction being liable to take place
between a pair of particles if they come to lie close enough to one
another. The density of particles at the initial time was chosen so
that the dynamics occur in a regime of mean free path, wherein a
typical particle meets a bounded number of other particles in a unit
of time. Speaking in rough terms, this choice of scaling causes the
effects of diffusion and interaction on the macroscopic evolution of
the system to be comparable. In common with much of non-equilibrium
statistical mechanics, we interpret the macroscopic behaviour of the
system in terms of the evolution of a small number of thermodynamic
parameters, in this case, the density of particles of a given mass, as
a function of macroscopic space and time.
In \cite{HR}, we  proved
that, when the initial number of particles is chosen to be high, this
density typically evolves as the solution of the Smoluchowski system
of PDE,
\setcounter{equation}{0}
\begin{equation}\label{syspde}
\frac{\partial{f_n}}{\partial t}(x,t)  = d(n)\Delta f_n (x,t)  +
Q^n_1(f) (x,t) -
Q^n_2(f) (x,t) . \qquad \qquad n = 1,2, \ldots
\end{equation}
 The first term on the
right-hand-side of (1.1) corresponds to the diffusion among
particles of mass $n$, with $d(n)$ being one-half of the diffusion
rate of such particles.
The terms in (1.1) corresponding to the interaction of pairs 
of particles are
given by the gain term
\begin{equation}\label{gainterm}
Q^n_1(f) (x,t)  =  \frac 12\sum_{m=1}^n \macp (m,n-m) f_m(x,t)
f_{n-m}(x,t),
\end{equation}
and the loss term
\begin{equation}\label{lossterm}
    Q^n_2(f) =   f_n (x,t) \sum_{m=1}^{\infty} \macp (m,n) f_m (x,t).
\end{equation}
Here, the collection of constants $\macp:\mathbb{N}^2 \to (0,\infty)$
quantify the macroscopic propensity of mass at a pair of values to
combine.

We will be concerned with weak solutions of the system (\ref{syspde}), defined by the equality of the left- and right-hand sides of (\ref{syspde}) after multiplication by $J_n:\R^d \times [0,\infty) \to [0,\infty)$, integration in space-time, and integration by parts. The equality is demanded over all choices of sequences of compactly supported smooth functions $\big\{ J_n: n \in \N \big\}$
such that only finitely many terms in the sequence are not identically zero.

The arguments of \cite{HR} were valid in the case where the dimension $d$
of the system was assumed to be at least three. The question of the
behaviour of such a system in two dimensions is significantly
different, and it is this topic that we address in this paper. 
We now
turn to describe the model in more detail, after which, we will
discuss the ways in which the two-dimensional case differs from that
of higher dimensions.

We will be working with a collection of microscopic models, each model
carrying an index $N \in \mathbb{N}$, this being the total number of
particles present in the system at the initial time. Each of these $N$
particles is independently assigned a random integer mass and placed
at the initial time at a random location whose law depends on that
mass. More precisely, we will be describing the state of the system at
any given moment in time by a configuration, by which we mean a map
$q: I_q \to \mathbb{R}^2 \times \mathbb{N}$, whose domain $I_q$ is
some finite set of a countable index set $I$. That is, if $i \in I_q$
has $q(i) = (x_i,m_i)$, then the system currently contains a particle
of mass $m_i$ at $x_i \in \mathbb{R}^2$. 
To define the initial configuration,
we choose  a sequence of continuous functions $\big\{
h_n: \mathbb{R}^2 \to [0,\infty), \, n \in \mathbb{N} \big\}$ that
must satisfy some conditions that we will shortly specify. We set
$Z = \sum_{n=1}^{\infty}{\int_{\mathbb{R}^2}}{h_n} \in (0,\infty)$, and
choose $N$ points in $\mathbb{N} \times \mathbb{R}^2$ indepedently
according to a law whose density at $(x,n)$ is equal to $h_n(x)/Z$.
Selecting arbitrarily a set of $N$ symbols $\{i_j: j \in \{1,\ldots,
N\} \}$ from $I$, we define the initial configuration $q_0$ by insisting
that $q_0(i_j)$ is equal to the $j$-th of the randomly chosen members
of $\mathbb{N} \times \mathbb{R}^2$.

Each particle moves according to an independent Brownian motion whose
diffusion rate $2d(m)$ depends on its mass $m \in \mathbb{N}$. As we
will explain later, we require some conditions on the choice of the
function $d: \mathbb{N} \to \mathbb{R}$, although the restriction
imposed by these conditions is far from prohibiting the 
 physically reasonable choice
where $d$ is decreasing. Any pair of particles that approach to within
a certain range of interaction are liable to coagulate, at which time,
they disappear from the system, to be replaced by a particle whose
mass is equal to the sum of the colliding particles, and whose
location is at some point nearby the place where the collision took
place. This range of interaction is taken to be equal to a parameter
$\epsilon$, whose dependence on the total particle number $N$ must be
stipulated. We make the choice $N = \vert \log \epsilon \vert Z$. This
will ensure that a particle randomly chosen from those initally
present experiences an expected number of collisions
in a given unit of time that remains bounded away from zero and
$\infty$ as $N$ is taken to be high. The effects of motion
and reaction determine the macroscopic evolution of the system to
comparable extents in this scaling.

We now describe the mathematical details of these dynamics. 
Let $F: \{ \mathbb{R}^2 \times \mathbb{N} \}^I \to \R$ 
denote a smooth function, 
whose domain is given the product topology. 
The dynamics is such that the action on $F$ of the model's 
infinitesimal generator $\mathbb{L}$ is given by   
\begin{displaymath}
   (\mathbb{L} F)(q)  = \mathbb{A}_0 F (q) +  \mathbb{A}_C F (q),
\end{displaymath}
where the diffusion and collision operators are given by
\begin{equation}\label{diff}
\mathbb{A}_0 F (q) =
\sum_{i \in I_q}{d(m_i) \Delta_{x_i} F }
\end{equation}
and
\begin{eqnarray}
\mathbb{A}_C F (q) & = &\frac 12
   \sum_{i,j \in I_q}{\epsilon^{-2} \big\vert \log \epsilon
   \big\vert^{-1} 
   V \Big( \frac{ x_i - x_j
   }{\epsilon} \Big) \micp(m_i,m_j)} \label{dyna} \\
   & & \qquad \qquad
\bigg[ \frac{m_i}{m_i +
   m_j} F \big( S^1_{i,j} q \big) + \frac{m_j}{m_i +
   m_j} F \big( S^2_{i,j} q \big) - F(q) \bigg]. \nonumber
\end{eqnarray}
Note that:
\begin{itemize}
\item the function $V: \mathbb{R}^2 \to [0,\infty)$ is assumed to be
H\"older continuous  of
compact support, and  with $\int_{\mathbb{R}^2}{V(x) dx} = 1$.
\item we denote by $S^1_{i,j}q$ that configuration formed from $q$ by
removing the indices $i$ and $j$ from $I_q$, and adding a new index
from $I$ to which  $S^1_{i,j}q$ assigns the value $(x_i,m_i + m_j)$. The
configuration  $S^2_{i,j}q$ is defined in the same way, except that it
assigns the value  $(x_j,m_i + m_j)$ to the new index. The specifics
of the collision event then are that the new particle
appears in one of the locations of the two particles being removed,
with the choice being made randomly with weights proportional to the
mass of the two colliding particles.
\end{itemize}
We will denote by $\mathbb{P}_N$ the measure on functions from
$t \in [0,\infty)$ to the configurations determined by the process at
time $t$. Its expectation will be denoted $\mathbb{E}_N$.

The form of the collision term in (\ref{dyna}) differs from that used
in the case of higher dimensions, in that the factor of $\big\vert
\log \epsilon \big\vert^{-1}$ is absent in the latter case. To explain
why we make this change, we firstly recall the reason for the form of
the collision operator in the case when $d \geq 3$. 
Suppose that, for some such choice of the dimension, two particles
$(x_i,m_i)$ and $(x_j,m_j)$ have, at some time $t_0$, just become liable to interact, in
the sense that the difference $x_i - x_j$ has become of order
$\epsilon$. This state of affairs is liable to persist for a time of
order $\epsilon^2$, but not much longer: for $d \geq 3$ and $C$ a
large constant, the Brownian
displacement $x_i - x_j$ would return a distance of $\epsilon$ from
the origin with only a small probability after a time of $C \epsilon^2$
after the moment $t_0$. 
This means
that, by choosing a form of collision dynamics in which the factor of
$\vert \log \epsilon \vert^{-1}$ is absent from (\ref{dyna}), we
ensure that the integral
\begin{displaymath}
   I_T = \int_{t =t_0}^{T}{\micp(m_i,m_j) \epsilon^{-2}
V \Big( \frac{x_i - x_j}{\epsilon} \Big)} dt,
\end{displaymath}
reaches its eventual value after a time of order $\epsilon^2$. Other
particles are unlikely to interfere with this pair in such a short
period of time, and, as such, we may neglect their influence. The
probability of collision between the pair before time $T$ is equal to
$1 - \exp \big\{ - I_T \big\}$. Thus, for $d \geq 3$, our choice of
dynamics is such that, among all the pairs of particles that at some
moment lie within $\epsilon$ of each other, the fraction that
eventually coagulate 
is bounded in $N$ away from $0$ and $1$, with this fraction
being close to $0$ or $1$ depending on whether the relevant constant
$\micp(m_i,m_j)$ is high or low.  

Turning to the planar case, note firstly that, 
in order that the probability of pair collision may remain
equal to $1 - \exp \big\{ - I_T \big\}$, we alter the definition of
$I_T$ by introducing a factor of $\big\vert \log \epsilon \big\vert^{-1}$.  
The two-dimensional case differs, because a planar Brownian motion
returns almost surely to any open set at indefinitely later times. As
such, the difference $x_i - x_j$ will endlessly re-enter the $\epsilon$-ball
centred at the origin, ensuring that $I_T \to \infty$ as $T \to
\infty$. In a system of two particles, their coagulation is
inevitable. In the system that we consider, where the regime of
constant mean free path has been selected, a pair of particles at
$\epsilon$ distance may find that their ongoing efforts to coagulate,
as measured by the increase of $I_T$, are interrupted by the arrival
of a third particle, the probability of appearance of such an intruder becoming
appreciable at a small time, independent of $\epsilon$, after that
at which the pair in question first came close to each other.   
The factor of  $\big\vert \log \epsilon \big\vert^{-1}$ that appears in
(\ref{dyna}) ensures that, during this short fixed time, the
probability of coagulation between the pair is of unit order, with the
constant $\micp(m_i,m_j)$ determining whether this probability is
high or low. To write a statement analogous to that for the higher
dimensional case: among the set of pairs of particles that are at some
moment at a distance of order $\epsilon$, the fraction that combine
with each other, rather than with some other particles, is bounded in
$N$ away from $0$ and $1$, with
the value of the constant
$\micp$ determining whether this fraction is high or low, similarly
to the earlier case. 

Our main result is conveniently expressed in terms of the empirical
measures on the locations of particles of a given mass. For each $n
\in \mathbb{N}$ and $t\in [0,\infty)$, we write
$g_n(dx,t)$ for the measure on $\mathbb{R}^2 $ given by
$$
g_n(dx,t)
=|\log\e|^{-1}\sum_{i\in I_{q(t)}}\d_{x_i(t)}(dx)1\!\!1\big(m_i(t)=n\big).$$

We write $g$ for the random measure on space-mass-time $\mathbb{R}^2 \times \N \times [0,\infty)$ such that, for each $t \geq 0$, its time-$t$ marginal $g(\cdot,t)$ is given by
$$
g(\cdot,t)
=  |\log\e|^{-1} \sum_{i\in I_{q(t)}}\d_{\big(x_i(t),m_i(t)\big)}.
$$

We also require a mild hypothesis on the diffusion coefficients
$d:\mathbb{N} \to (0,\infty)$ (see the first remark after Theorem 1.1 below). Namely, we
suppose that there exists a function $\g : \mathbb{N}^2 \to 
(0,\infty)$ such that $\micp \le \g$,
with $\g$ satisfying
\begin{equation}\label{hypo}
n_2 \g \big( n_1, n_2 + n_3 \big) \max{ \Big\{ 1, \Big[ \frac{d (
n_2 + n_3 )}{d(n_2)}
\Big]^{3}\Big\}} \leq \big( n_2 + n_3 \big) \g(n_1,n_2).
\end{equation}

The initial random configuration of $N$ particles is formed by
scattering particles of numerous masses independently in
$\mathbb{R}^2$ according to densities that are prescribed for each
mass. These densities will be chosen as continuous functions
$\{ h_n: \mathbb{R}^2 \to [0,\infty), n \in \mathbb{N} \}$, and should
satisfy some fairly weak bounds. To be specific, we insist that
\begin{itemize}
\item $k\in L^1(\R^2)$ and $\bar k\in L^{\i}_{loc}({\mathbb{R}^2})$ where
$k:=\sum_{n=1}^{\i} nh_n$ and 
$$
\bar k(x)=\iint_{|x_1-x_2-y|\le 1}\big|\log |x_1-x_2-y|\big|k(x_1)k(x_2)dx_1dx_2.
$$
\item For every $m$,  $\sum_{n=1}^{\infty}d(n)^{2/3}\g(n,m){\hat h_n} \in
L^\infty_{loc}({\mathbb{R}^2})$ for
$\hat h_n(x)=\int  h_n(y)|x-y|^{-4/3}dy$.
\item For every $m$, $\sum_{n=1}^{\infty}d(n)^{3/4}\g(n,m){\tilde h_n} \in
L^\infty_{loc}({\mathbb{R}^2})$ for
$\tilde h_n(x)=\int  h_n(y)|x-y|^{-3/2}dy$.
\end{itemize}
We then set
$Z=\sum_{n=1}^{\infty}\int_{\mathbb{R}^2}{h_n} \in (0,\infty)$ and choose $N$
points
in $\mathbb{N} \times \mathbb{R}^2$ independently according to a law
whose density at $(x,n)$ is equal to $h_n(x)/Z$.
Selecting arbitrarily a set of $N$ symbols $\{i_j: j \in \{1,\ldots,
N\} \}$ from $I$, we define the initial configuration $q_0$ by insisting
that $q_0(i_j)$ is equal to the $j$-th of the randomly chosen members
of $\mathbb{N} \times \mathbb{R}^2$.
\newline
\noindent
{\bf Remark.}
It is not hard to show that our assumptions on the initial data $\{h_n\}$
are satisfied if $k$ is bounded, has a bounded support,
$d(\cdot)$ is bounded and $\gamma(n,m)\le C(m)n$ for a function $C(\cdot)$.
 Indeed if $k$ is bounded and has a bounded support, then
$\bar k L^\i_{loc}$ and $\hat h_n,\tilde h_n\in 
L^\i_{loc}$ for every
$n$.  It is worth mentioning that if $k$
belongs to the negative Sobolev Space $H^{-1}=W^{-1,2}$, then $\bar k\in L^\i$.

The main theorem is now stated.
\begin{theorem}\label{thmo} 
 Let
$\mathcal P_N$ denote the law on measures on $\R^2 \times \N \times [0,\infty)$ given by the law of $g$ under $\mathbb P_N$; 
recall that $\epsilon$ is related to $N$ by means of the
formula $N
|\log\e|^{-1} = Z$, with the constant $Z \in (0,\infty)$ being
given by the expression $Z = \sum_{n \in
\mathbb{N}}{\int_{\mathbb{R}^2}h_n}$.

The sequence $\{\mathcal P_N\}$ is 
tight. 
 Moreover, any limit point $\mathcal P$ of the sequence $\{\mathcal P_N\}$ is concentrated on the space of measures taking the form $\sum_{n=0}^{\infty} f_n(x,t) dx \times \delta_n \times dt$ where 
$\{f_n: n \in \mathbb N\}$ ranges over weak solutions of (\ref{syspde}) that satisfy
the initial condition $f_n(\cdot,0) = h_n(\cdot)$; recall that  the collection of
constants $\macp: \mathbb{N}^2 \to [0,\infty)$ is given by 
\begin{equation}\label{recp}
 {\b}(n,m)= \frac{2\pi \big( d(n) + d(m) \big)\a(n,m)}{2\pi \big( d(n)
+ d(m) \big) + \micp(n,m)}.
\end{equation}
\end{theorem}
Note that convergence in Theorem~\ref{thmo} is asserted only subsequentially and to a limiting object which may be a random superposition of weak solutions of~(\ref{syspde}). The need for such a weak statement of convergence disappears in the case that uniqueness of the weak solution of~(\ref{syspde}) are known.
Some such conditions are provided by Proposition 2.6 of \cite{wrzd}. Shortly after the present paper originally appeared, we proved uniqueness in a reasonably general setting. The next result is a consequence of the main theorems of~\cite{momentbounds} as explained in Remark~1.2 of that paper.
\begin{prop}\label{propuniqueness}
Let the dimension satisfy $d \geq 1$. 
For $a,b > 0$ such that $a + b < 1$, and for positive constants $c_1$ and $c_2$, assume that 
$\macp(n,m)\le c_1(n^a+m^a)$ and $d(n)\ge c_2 n^{-b}$ for
all $n,m\in \N$. Also assume that $d:\N \to (0,\infty)$ is non-increasing. There exists $e > 0$ such that $\sum_nn^e \|h_n\|_{L^\infty(\R^d)}<\i$
 and $\|\sum_nn^e h_n\|_{L^1(\R^d)}<\i$ imply that~(\ref{syspde}) has a unique weak solution. (In fact,  this solution conserves mass, in the sense that  $I:[0,\infty) \to [0,\infty)$ given by $I(t)  = \sum_{m \in \mathbb{N}} m \int_{\mathbb{R}^d} f_m(x,t)dx$ satisfies $I(t) = I(0)$ for all $t \in [0,\infty)$.)
\end{prop}

Theorem~\ref{thmo} and Proposition~\ref{propuniqueness} permit convergence of the empirical measures to be asserted in a more satisfying sense:
\begin{cor}\label{coro}
Suppose that the assumptions of 
Theorem~\ref{thmo} and Proposition~\ref{propuniqueness} are in force.

Let $J: \mathbb{R}^2 
\times [0,\infty) \to \mathbb{R}$ be a 
 bounded and continuous test function. Then, for each
$n \in \mathbb{N}$ and $T \in (0,\infty)$,
\begin{equation}\label{resu}
\limsup_{N \to \infty } \mathbb{E}_N  \bigg\vert \int_{[0,T]} \int_{\mathbb{R}^2}
J(x,t) \big( g_n(dx,t) - f_n (x,t) dx \big) dt \bigg\vert = \ 0,
\end{equation}
where again
 $N
|\log\e|^{-1} = Z$, with  $Z = \sum_{n \in
\mathbb{N}}{\int_{\mathbb{R}^2}h_n}$. In (\ref{resu}),
    $\{ f_n: \mathbb{R}^2 \times [0,\infty) \to [0,\infty), n \in \mathbb{N} \}$
denotes the unique weak solution to the system of
partial differential equations (\ref{syspde}) with the initial data $f_n(\cdot,0)
=h_n(\cdot)$.
\end{cor}
{\bf Remarks} 
\begin{itemize}
\item
Included in the space of parameter values that satisfy
(\ref{hypo}) is the case where the diffusion rate $d$ is a decreasing
function of the mass, and the coagulation propensities $\micp$
satisfy $\micp \big( n , m \big) \leq C nm$. In fact for a 
nonincreasing $d(\cdot)$, the condition (\ref{hypo})
is equivalent to saying that $\micp(n,m)\le C(n)m$ for a function $C(n)$.
Also, if the microscopic coagulation
rate $\micp$ is identically constant, then the condition (\ref{hypo}) 
is equivalent to saying that the function
$d(n)n^{-1/4}$ is nonincreasing.
\item Note that the macroscopic coagulation propensities $\b$  depend only
on the total integral of $V$ that is assumed to be $1$ for convenience. 
However when the dimension is $3$ or more the propensity $\b(n,m)$ does depend on $V$ in
a nontrivial way and
is given as $\a(n,m)\int (1+u)Vdx$, where $u$ solves the PDE $\Delta u=\tau (1+u)V$ with
$\tau=\a(n,m)/(d(n)+d(m))$.   
\item Our technique of proof also yields a kinetic limit derivation
for the model in which particles are assumed to have a range of
interaction that is mass-dependent. To give an example of such a
variant, suppose that each
particle of mass $m$ has a radius $r(m)$, where $r(m) =\sqrt m$.
 We stipulate that particles of mass $m$ and
$n$ are liable to react when their displacement reaches the order of
$\big( r(m) + r(n) \big) \epsilon$. More precisely, we modify the
definition (\ref{dyna}) of the collision operator $\mathbb{A}_{C}$ by
replacing the appearance of $V $ by $ \big( r(n) + r(m) \big)^{-2}  V \big(
\cdot/(r(n) + r(m))\big)$, (the factor that multiplies $V$ being
introduced so that, roughly speaking, the altered collision mechanism
respects the spatial-temporal scaling of Brownian motion).  Theorem 1.1 is still valid 
for this modified model with the same macroscopic coagulation propensities $\b$. 
This is  in sharp contrast with the case $d\ge 3$ for which the mass dependence 
affects the macroscopic coagulation propensities~$\b$.
\end{itemize}

In common with the proof for $d \geq 3$, a central element in deriving
Theorem \ref{thmo} is establishing that, at any given moment after the
initial time, the presence of a particle of some given mass at some
fixed point in space significantly affects the likelihood of a
particle being at some other point in space only if that other point is at a
short distance from the first particle. That is, on distances of short
order, the presence of a particle makes it less likely to find another
nearby, because the pair would have been liable to coagulate shortly
beforehand. However, the distribution of particles at a given time is
similar to one in which they were scattered independently, except for
this short-range repulsion. The following proposition, whose form differs from that in the case $d \geq 3$ only in its scaling factor, formalises this assertion.
\begin{prop}\label{szprop}
Set 
\begin{equation}\label{qudef}
Q = |\log\e|^{-2} \sum_{(i,j) \in
I_{q}}{\micp(m_i,m_j) V_{\epsilon}(x_i - x_j) J(x_i,m_i,t)
    \overline{J}(x_j,m_j,t)},
\end{equation}
where $J,\overline{J}: \mathbb{R}^2 \times \mathbb{N}
\times [0,\infty) \to [0,\infty)$ are test functions satisfying the
same conditions as those stated in Theorem \ref{thmo}.
We also assume that
$J(x,m,t)=0$ unless $m=M_1$ and $\overline J(x,m,t)=0$ unless $m=M_2$.
Let $\eta:
\mathbb{R}^2 \to [0,\infty)$ denote a smooth function of compact
support for which $\int_{\mathbb{R}^2}{\eta(x)dx} = 1$.  
We have that
\begin{eqnarray}\label{szrec}
    & & \int_0^T  Q (t) dt \\
    & = &
\int_0^T{} dt \int_{\mathbb{R}^2}
d\omega  \, \b(M_1,M_2)  J(\omega,M_1,t)
\overline{J}(\omega,M_2,t)  \bigg[ |\log \e|^{-1} \sum_{i \in I_q: m_i = M_1}{
\delta^{-2} \eta \Big( \frac{x_i - \omega}{\delta} \Big) }
\bigg] \nonumber \\
& & \qquad \quad \bigg[ |\log\e|^{-1} \sum_{j \in I_q; m_j = M_2}{
\delta^{-2} \eta \Big(\frac{x_j - \omega}{\delta} \Big)
} \bigg]  \, + \rm{Err} \big( \epsilon, \delta
\big), \nonumber
\end{eqnarray}
where the constants ${\b}: \mathbb{N}^2 \to [0,\infty)$ were
defined in (\ref{recp}), and
where the function $\rm{Err}$ satisfies
$$
  \lim_{\delta \downarrow 0}\limsup_{\epsilon\downarrow 0} \mathbb{E}_N
  \big\vert Err \big( \epsilon, \delta \big) \big\vert = 0.
$$
\end{prop} 
Why is this statement a mathematical rendering of the claim discussed
before it was made? 
The quantity $\int_0^T Q(t) dt$ can be thought of as the total
propensity of particles to combine during 
the interval of time $[0,T]$. Proposition~\ref{szprop} asserts it may
be approximated by a time-averaged product of empiricial
approximations to the density of particles (of the appropriate
mass). That is, particles are arranged independently enough near most
of the collision events that the rate of these collisions is roughly
proportional to that arising in a system in which particles are scattered
indepedently at random according to densities given by measuring the
system in question on scale $\delta$ that is much larger than the
reaction range $\epsilon$. There is, however, a constant of proportion
corresponding to the change from microscopic reaction propensity
$\micp$ appearing in the definition of $Q$ to its macroscopic
counterpart ${\b}$. Its presence may be explained by
the negative short-range correlation between particles discussed
before the statement of Proposition \ref{szprop}.
 
The analogue of Theorem \ref{thmo} that appears
in \cite{HR} for the case $d \geq 3$ is derived as a consequence of
Proposition \ref{szprop}. In Section 2 of \cite{HR}, a sketch of the proof of
Theorem \ref{thmo} may be found. The details of the derivation of
Theorem \ref{thmo} from Proposition \ref{szprop} do not differ in the
two-dimensional case, so that we do not present these arguments again
in this paper. Our task here is rather to present a detailed
derivation of Proposition \ref{szprop} in the case when $d=2$. Before
reading further, however, the reader may wish to consult Section 2 of \cite{HR}. We refer the reader to \cite{HR} also for a discussion of previous work related to the problem. Here, we mention only Sznitman
\cite{sznitman}, in which a model of Brownian spheres that annihilate as soon as
they touch is studied. The partial differential equation by which the density of
particles evolves was derived for the kinetic limit, in each dimension
$d \geq 2$. In this work, the macroscopic annihilation rate is exactly ${2\pi}$
when the dimension is $2$. This is compatible with our main results because if $d(\cdot)$
is identically $1/2$ and $\a\to\i$, the macroscopic coagulation rate $\b$
approches $2\pi$. Note that our model approximates the {\it hard core} model 
as $\a$ gets large.

\medskip

\noindent{\bf Acknowledgment.} The first author would like to express his thanks to James Norris for introducing him to the topic of diffusive coagulating systems and for valuable discussions.
\end{section}
\begin{section}{Establishing the Stosszahlansatz}\label{sectwo}
\setcounter{equation}{0}

For any given pair $(n,m) \in \mathbb{N}^2$,  $\ue =
\ue_{n,m}: \mathbb{R}^2 \to [0,\infty)$ is a
function whose existence is ensured by Theorem \ref{th3.1} that lies in $C^2
\big( \mathbb{R}^2 \big)$ satisfying
\begin{equation}\label{eq3.1}
     \ue_{n,m} \big( x \big)  = \frac 1{2\pi}\tau \big( n,m \big) \int
     \log |x-y|\Big[
V_{\epsilon} \big(
   y \big) \ue_{n,m} \big( y \big) +  V^{\epsilon} \big( y \big) \Big] {\rm d} y,
\end{equation}
where ${\tau}(n,m)=\a(n,m)/(d(n)+d(m))$. As a consequence,
\begin{displaymath}
\big(d(n)+d(m)\big)\Delta \ue_{n,m} \big( x \big)  = \a \big( n,m \big) \Big[
V_{\epsilon} \big(
    x \big) \ue_{n,m} \big( x \big) +  V^{\epsilon} \big( x \big) \Big].   
\end{displaymath}
We are using the notations
\begin{displaymath}
    V^{\epsilon}(x)  = \epsilon^{-2} V \Big( \frac{x}{\epsilon} \Big)
\end{displaymath}
and
\begin{displaymath}
    V_{\epsilon}(x)  = \epsilon^{-2}|\log\eps|^{-1} V \Big( 
\frac{x}{\epsilon} \Big).
\end{displaymath}

We present the conditions on the two test functions $J,\overline{J}:
\mathbb{R}^2 \times \mathbb{N} \times [0,\i)
\to \R$ that appear in Proposition~\ref{szprop}. It suffices to work with
functions that take
non-zero values for only one value in the second argument, such
functions measuring the presence of particles of a given mass. By a temporary
abuse of notation, we write
\begin{eqnarray}
J(x,M_1,t) & = & J(x,t) 1 \! \! 1  \big\{ m = M_1 \big\} \nonumber \\
\overline{J}(x,M_1,t) & = & \overline{J}(x,t) 1 \! \! 1  \big\{ m =
M_2 \big\}, \nonumber
\end{eqnarray}
where on the right-hand-side, $J$ and $\overline{J}$ denote smooth
maps from $\mathbb{R}^2 \times [0,\infty)$ to $\R$ of compact support. We will
suppress the appearance of the $t$-variable when writing the arguments
of $J$ and $\overline{J}$.

In seeking to verify the Stosszahlansatz, we define
\begin{equation}\label{exzo}
X_z (q) = |\log\epsilon|^{-2} \sum_{i,j \in I_q}{
\ue_{M_1,M_2} (x_i- x_j + z) J (x_i,M_1,t) \overline{J}(x_j,M_2,t)} 
 1\! \! 1 \big\{ m_i = M_1,m_j = M_2 \big\} .
\end{equation}
The relevance of the expression (2.2) for our purposes 
is that the term $Q$ and its variations appear as we apply the infinitesimal generaor on 
the expression $X_z-X_0$. We refer the reader to Section 2 of \cite{HR} for some
 heuristic justification of the special form of $X_z$. 
 
Numerous terms arise when the operators $\mathbb{A}_0$ and
$\mathbb{A}_C$ act on the expression $X_z - X_0$ (recall that the
functions of configurations $X_z$, indexed by $z \in \mathbb{R}^2$,
were defined in (2.2)). We now label these terms.
Unless stated otherwise, we will adopt a notation whereby all the
index labels appearing in sums should be taken to be distinct. This
includes the case of multiple sums. For example, $\sum_{k,l \in I_q}
\sum_{i \in I_q} f(x_k,x_l,x_i)$ denotes the sum of the evaluation of the
function $f$ over all arguments that are triples $(x_k,x_l,x_i)$ where
$k$, $l$ and $i$ are distinct indices in $I$. Note also that, unless
otherwise stated, whenever
the symbol $\ue$ appears in a summand, we mean $\ue_{M_1,M_2}$.

Firstly, we label those terms arising
from the action of the diffusion operator. To do so, note that, for a
time-dependent functional $F$ of the configuration space, this action
is given by
\begin{displaymath}
    \Big( \frac{\partial}{\partial t} +  \mathbb{A}_0 \Big) F =
\frac{\partial}{\partial t} F + \sum_{i \in I_q}{d \big(m_i)
    \Delta_{x_i} } F.
\end{displaymath}
Thus, we label as follows:
$$
\left(\frac{\partial}{\partial t}+\mathbb{A}_0\right) (X_z - X_0)( q(t) ) =
H_{11} + H_{12} + H_{13} + H_{14}
+ H_2 + H_3 + H_4,
$$
with
\begin{eqnarray}
H_{11} & = & \sparsize \sum_{i,j \in I_q}{\alpha(m_i,m_j)
    \Big[ V^{\epsilon} \big( x_i - x_j + z  \big) -
V^{\epsilon} \big(   x_i - x_j  \big) \Big]} J(x_i,m_i,t)
    \overline{J}(x_j,m_j,t) \nonumber \\
H_{12} & = & - \sparsize \sum_{i,j \in I_q}{\alpha(m_i,m_j)
    V_{\epsilon} \big( x_i - x_j  \big) \ue \big( x_i - x_j
    \big)  J(x_i,m_i,t) \overline{J}(x_j,m_j,t) } \nonumber \\
H_{13} & = &  \sparsize \sum_{i,j \in I_q}{\alpha(m_i,m_j)
     V_{\epsilon} \big( x_i - x_j + z \big) \ue \big( x_i - x_j+z
    \big)  J(x_i,m_i,t) \overline{J}(x_j,m_j,t) }, \nonumber \\
H_{14} & = &  \sparsize \sum_{i,j \in I_q}{
     \Big[ \ue \big( x_i - x_j + z \big) -  \ue \big(
    x_i - x_j \big) \Big]} \nonumber \\
    & & \qquad \quad \Big[    J_t \big(x_i,m_i,t\big)
    \overline{J}(x_j,m_j,t) +   J \big(x_i,m_i,t\big)
      \overline{J}_t(x_j,m_j,t) \Big] , \nonumber
\end{eqnarray}
along with
\begin{eqnarray}
H_2 & = & 2 \sparsize \sum_{i,j \in I_q}{ d(m_i)
      \overline{J}(x_j,m_j,t)} \nonumber \\
    & & \qquad \qquad \qquad \quad \Big[  \ue_x (x_i - x_j + z ) - \ue_x
(x_i - x_j ) \Big]
    \cdot J_x(x_i,m_i,t), \nonumber \\
H_3 & = & - 2 \sparsize \sum_{i,j \in I_q}{ d(m_j)
      J(x_i,m_i,t)} \nonumber \\
     & & \qquad \qquad \qquad \quad  \Big[  \ue_x (x_i - x_j + z ) -
     \ue_x (x_i - x_j ) \Big]
    \cdot \overline{J}_x (x_j,m_j,t) , \nonumber
\end{eqnarray}
and
\begin{eqnarray}
H_4 & = & \sparsize \sum_{i,j \in I_q}{} \Big[
    \ue (x_i - x_j + z ) - \ue (x_i - x_j ) \Big]
\nonumber\\
    & &  \qquad \qquad \quad  \Big[ d(m_i) \Delta_x J(x_i,m_i,t)
    \overline{J}(x_j,m_j,t) +  d(m_j) J(x_i,m_i,t)
    \Delta_x \overline{J}(x_j,m_j,t) \Big] , \nonumber
\end{eqnarray}
where $f_x$ denotes the gradient of $f$, and $\cdot$ the scalar
product. As for those terms arising from
the action of the collision operator,
\begin{displaymath}
\mathbb{A}_C (X_z - X_0)(q) = G_z (1) + G_z (2) - G_0 (1) - G_0 (2),
\end{displaymath}
where $G_z (1)$ is set equal to
\begin{eqnarray}
& &\frac 12 \sum_{k,l \in
I_q}{\alpha(m_k,m_l) V_{\epsilon}( x_k - x_l)} \sparsize \sum_{i
\in I_q}{} \nonumber \\
    & & \quad \bigg\{
    \frac{m_k}{m_k + m_l} \Big[  \ue ( x_k - x_i
    + z) J( x_k , m_k + m_l,t)
    \overline{J}(x_i,m_i,t) \nonumber \\
    & & \qquad \qquad \qquad  + \, \ue ( x_i - x_k
    + z) J( x_i , m_i ,t) \overline{J}(x_k,m_k + m_l,t)  \Big]
     \nonumber \\
    & & \quad \, \, + \, \frac{m_l}{m_k + m_l}  \Big[ \ue ( x_l - x_i
    + z) J( x_l , m_k + m_l,t)
    \overline{J}(x_i,m_i,t) \nonumber \\
    & &  \qquad \qquad \qquad  + \, \ue ( x_i - x_l
    + z) J( x_i , m_i ,t) \overline{J}(x_l,m_k + m_l,t) \Big]
     \nonumber \\
    & & \qquad  - \ \
    \Big[  \ue (x_k - x_i + z) J( x_k , m_k,t)
    \overline{J}(x_i,m_i,t) \nonumber \\
    & &  \qquad \qquad \qquad \qquad + \,  \ue (x_i - x_k + z)
    J( x_i , m_i ,t) \overline{J}(x_k,m_k,t) \Big]
     \nonumber \\
    & & \qquad - \ \
     \Big[  \ue (x_l - x_i + z) J( x_l , m_l,t)
    \overline{J}(x_i,m_i,t) \nonumber \\
    & & \qquad \qquad \qquad \qquad + \,  \ue (x_i - x_l + z) J( x_i ,
    m_i ,t) \overline{J}(x_l,m_l,t) \Big] \bigg\}
     \nonumber ,
\end{eqnarray}
and where
\begin{eqnarray}\label{eqnzeta}
G_z (2) & = & - \sparsize \sum_{k,l \in
I_q}{ \alpha(m_k,m_l)  V_{\epsilon}(x_k - x_l)} \\
    & & \qquad \quad
\ue (x_k -
x_l + z ) J(x_k,m_k,t)  \overline{J}(x_l,m_l,t). \nonumber
\end{eqnarray}
The terms in $G_z(1)$
arise from the changes in the functional $X_z$ when a collision
occurs due to the influence of the appearance and disppearance of
particles on other particles that are not directly involved. Those in
$G_z(2)$ are due to the absence after collision of the summand in
$X_z$ indexed by the colliding particles.

Note that
\begin{equation}\label{igt}
    H_{12} + G_0 \big( 2 \big) = 0.
\end{equation}
The process $\big\{ \big( X_z - X_0 \big)(t): t \geq 0 \big\}$ satisfies
\begin{eqnarray}
    \big( X_z - X_0 \big) \big( T \big) & = &
    \big( X_z - X_0 \big) \big( 0 \big) +
\int_{0}^{T}{ \Big( \frac{\partial}{\partial t} + \mathbb{A}_0 \Big)
\big(X_z - X_0\big) (t)dt} \label{hact} \\
    & & \qquad + \,
\int_{0}^{T}{ \mathbb{A}_C (X_z - X_0) (t)dt} \, + \, M(T), \nonumber
\end{eqnarray}
with $\big\{ M(t) : t \geq 0 \big\}$ being a martingale.
By using the labels for the various terms
that we just introduced, we find from (\ref{hact}) by use of (\ref{igt}) that
\begin{eqnarray}
& & \Big\vert
\int_{0}^{T}{ H_{11} \big(t\big) dt} +
\int_{0}^{T}{ H_{13} \big(t\big) dt} \Big\vert \nonumber \\
& \leq &  \big\vert X_z - X_0 \big\vert \big( q(T) \big) + \big\vert
X_z - X_0 \big\vert \big( q(0) \big) \nonumber \\
& & + \int_{0}^{T}{\big\vert H_{14}  \big\vert (t) 
dt}+\int_{0}^{T}{\big\vert H_2
  \big\vert (t) dt} +
\int_{0}^{T}{\big\vert  H_3  \big\vert (t) dt}  +
\int_{0}^{T}{\big\vert  H_4 \big\vert (t) dt} \nonumber \\
& & +
\int_{0}^{T}{  \big\vert G_z(1) - G_0(1) \big\vert (t) dt} +
\int_{0}^{T}{  \big\vert G_z(2) \big\vert (t) dt} \, + \, \big\vert
M(T) \big\vert. \label{abc}
\end{eqnarray}
Since $J$ is of compact support, we have that $X_z(q(T))=0$ for
$T$ sufficiently large.
We aim to prove the following estimates: for each $T > 0$,
\begin{eqnarray}
& & \int_{0}^{T}{ \mathbb{E}_N \vert H_{14} \vert (t) dt} \leq C
|z|^{1/2}\big|\log |z|\big|, \nonumber \\
& & \int_{0}^{T}{ \mathbb{E}_N \vert H_2 \vert (t) dt} \leq C \vert z
\vert^{1/9}\big|\log|z|\big|,
\nonumber \\
& & \int_{0}^{T}{ \mathbb{E}_N \vert H_3 \vert (t) dt} \leq C \vert z
\vert^{1/9}\big|\log|z|\big|, \nonumber \\
& & \int_{0}^{T}{ \mathbb{E}_N \vert H_4 \vert (t) dt} \leq  C
|z|^{1/2}\big|\log |z|\big|, \nonumber \\
& & \int_{0}^{T}{ \mathbb{E}_N \vert G_z(1) - G_0 (1) \vert (t) dt} \leq
C  |z|^{1/2}\big|\log |z|\big| | , \nonumber \\
& & \int_{0}^{T}{ \mathbb{E}_N \vert G_z (2) \vert (t) dt} \leq C 
\big\vert \log \vert z
\vert \big\vert  |\log\e|^{-1},
\nonumber\\
& &\mathbb{E}_N|X_z-X_0(0)|\le C|z|. \label{jayeye}
\end{eqnarray}
Later, we apply the limit $\vert z \vert \to 0$  after sending $\e$ to $0$.
We will also show that, for
each $T \in (0,\infty)$,
\begin{equation}\label{martest}
\mathbb{E}_N \Big[ M \big( T \big)^2 \Big] \leq
  C \big\vert \log \eps \big\vert^{-1} .
\end{equation}
\subsection{Lemmas bounding collision propensity}\label{stone}
In this subsection, we discuss three lemmas that in essence serve as
the backbone of the proof of the various inequalities that appear in (2.7).
These lemmas allow us to reduce the proof to a calculation involving the
initial configuarions for which the independence of particles and our assumptions on the 
initial densities can be used. In fact the proof of Lemmas 2.1 and 2.3 is 
very similar to the corresponding Lemmas 3.1 and 3.3 of \cite{HR}. For this reason,
their proofs are omitted. It is Lemma 2.2 that is somewhat different from what we have in 
\cite{HR} as Lemma 3.2 and we provide a detailed proof for it. In fact this difference explains 
to some extent a major technical difficulty that is two dimensional and is not encountered 
when the dimension is 3 or more. To explain this further, let us observe that if 
the dimension $d$ is 3 or more and
$J$ is a nonnegative function, then we can find a solution to the Poisson equation
$-\Delta H=J$ that satisfies $H\ge 0$. Indeed the solution $H$ is defined by
$$
c_0(d)\int |x-y|^{2-d}J(y)dy,
$$
where $c_0(d)=(d(d-2)\o(d))^{-1}$ with $\o(d)$ denoting the volume of the unit ball in 
$\R^d$. This is no longer true in dimension 2 because the solution is given by
$$
-\frac 1{2\pi}\int \log|x-y|\  J(y)dy.
$$
This causes some difficulty in treating various terms that appear in (2.7). To get around this,
let us define
\begin{equation}\label{eq2.8}
H(x)=-\frac 1{2\pi}\int_{|x-y|\le 1} \log|x-y|\  J(y)dy.
\end{equation}
We now have that $-\Delta H=J-\tilde J$ where 
\begin{equation}\label{eq2.9}
  \tilde{J}(x) = \frac 1{2\pi}\int_{\vert z \vert = 1}{J(x-z) dS(z)},
\end{equation}
where $dS$ denotes the $1$-Lebesgue measure on the unit circle $S^1$.
The point is that by using Lemma 2.2, we reduce bounding an expression involving
$J$ to an expression involving $H$ at time $t=0$, and a similar expression involving
$\tilde J$. Since the funcion $\tilde J$ is an average of $J$, we have an easier 
task to bound the expression involving $\tilde J$. In the case of the terms $H_2$ and $H_3$,
we need to apply this process three times so that the final $\tilde J$ has a simple
pointwise bound.
Our three lemmas are:

\begin{lemma}\label{lembc}
For any $T \in [0,\infty)$,
$$
\parsize \mathbb{E}_N \int_0^T dt \sum_{i,j \in I_q}{\alpha(m_i,m_j)}
    V_{\epsilon} \big( x_i - x_j \big) \leq 2Z.
$$
\end{lemma}
\begin{lemma}\label{lemtwop}
Let $J: \mathbb{R}^2 \to [0,\infty)$ be continuous, and let $H:
\mathbb{R}^2 \to [0,\infty)$ be given by (2.9)
We also define $\tilde{J}: \mathbb{R}^2 \to \mathbb{R}$ according to
(2.10).
Then we have the following inequality,
\begin{eqnarray}
& &    \big\vert \log \epsilon \big\vert^{-2} \mathbb{E}_N \int_0^T
\sum_{i,j \in I_q}{J(x_i - x_j)
m_i m_j (d(m_i) + d(m_j))} dt  \nonumber \\
& &  +  \big\vert \log \epsilon \big\vert^{-2} \mathbb{E}_N \int_0^T
\sum_{i,j \in I_q}{V_{\e}(x_i-x_j)\alpha(m_i,m_j)H(x_i - x_j)
m_i m_j } dt  \nonumber \\
\ \ \ & \leq &  \big\vert \log \epsilon \big\vert^{-2}
  \mathbb{E}_N \sum_{i,j \in I_{q(0)}}{H(x_i - x_j)
m_i m_j} \nonumber \\
\ \ \  & &  \, + \,
   \big\vert \log \epsilon \big\vert^{-2} \mathbb{E}_N \int_0^T
\sum_{i,j \in I_q}{\tilde{J}(x_i - x_j)
m_i m_j (d(m_i) + d(m_j))} dt. \nonumber
\end{eqnarray}
\end{lemma}
\begin{lemma}\label{lemthreep}
Assume that the function $\g:\mathbb{N}^2\to (0,\i)$ satisfies
\begin{equation}\label{hypo2}
n_2 \g \big( n_1, n_2 + n_3 \big) \max{ \Big\{ 1, \Big[ \frac{d (
n_2 + n_3 )}{d(n_2)}
\Big]^{2} \Big\}} \leq \big( n_2 + n_3 \big) \g(n_1,n_2),
\end{equation}
There exists a collection of constants $C: \mathbb{N}^2 \to (0,\infty)$,
such that, for any smooth function $J : \mathbb{R}^{4} \to [0,\infty)$,
and any given $n_1,n_3 \in \mathbb{N}$,
\begin{eqnarray}
& &
\mathbb{E}_N{\int_0^T{dt \sum_{k,l,i \in I_{q(t)}}{\g (m_i,m_j)
V_{\epsilon} \big(
    x_i - x_j \big) J ( x_i , x_k ) 1 \! \! 1 \{ m_i = n_1, m_k =
n_3 \} }}} \nonumber \\
    & \leq &  C_{n_1,n_3}\big\vert \log \eps \big\vert^{3}
    \sum_{n_2 \in \mathbb{N}}{\int{A^{\epsilon}_{n_1,n_2,n_3} \big( x_1,x_2,x_3
    \big)h_{n_1}(x_1)h_{n_2}(x_2)h_{n_3}(x_3) dx_1 dx_2 dx_3}},
\end{eqnarray}
where, also given $\epsilon > 0$ and $n_2 \in \mathbb{N}$, the function
$ A_{n_1,n_2,n_3}^{\epsilon}: \mathbb{R}^{6} \to
[0,\infty) $ is defined by
\begin{eqnarray}
    & & \big( d(n_1 )\Delta_{x_1} + d(n_2) \Delta_{x_2} + d(n_3)
\Delta_{x_3} \big)
    A_{n_1,n_2,n_3}^{\epsilon} (x_1,x_2,x_3) \nonumber \\
    & & \qquad = -   \g (n_1,n_2)
    V_{\epsilon} ( x_1 - x_2 ) J ( x_1 ,
    x_3 ) .
\end{eqnarray}
\end{lemma}

It is worth mentioning that the function $A_{n_1,n_2,n_3}^{\epsilon}(x_1,x_2,x_3)$
of Lemma 2.3 is given by
\[c_0 ( 6) \int_{\mathbb{R}^2}  \int_{\mathbb{R}^2}
\int_{\mathbb{R}^2}
     \bigg( \frac{\vert x_1 - z \vert^2}{d(n_1)}  + \frac{\vert x_2 -
y \vert^2}{d(n_2)} + \frac{\vert x_3 - y' \vert^2}{d(n_3)}
\bigg)^{-{2}} \g (n_1,n_2) J ( z , y' )  V_{\epsilon} ( z - y )
dz dy dy', 
\]
where $c_0(d)=(d(d-2)\omega(d))^{-1}$ with $\omega(d)$ denoting the volume of the
unit ball in $\R^d$. Note that for Lemma 2.3 we are dealing with a solution to
a Laplace type equation in $\R^6$ as opposed to Lemma 2.2 for which the pecularity
of the Laplace equation in $\R^2$ played a role. This is why the proof of \cite{HR}
in the case of Lemma 2.3 can be repeated line by line. 

{\bf Proof of Lemma \ref{lemtwop}}
Set
$$
X_q = \big\vert \log \epsilon \big\vert^{-2} \sum_{i,j \in I_q}{H(x_i 
- x_j)m_i m_j}.
$$
Recall the mechanism of the dynamics at collision:  the location of
the newly created particle is one of the
two locations of the colliding particles, with weights proportional to
the masses of the incident particles. We see that when $\mathbb{A}_C$ acts on
$X_q$, all those terms indexed by pairs of particles one of which is
not involved in the collision cancel. Thus,
\begin{equation}\label{locc}
\mathbb{A}_C X = - \big\vert \log \epsilon \big\vert^{-2} \sum_{i,j 
\in I_q}{ V_{\epsilon} ( x_i
- x_j ) \alpha(m_i,m_j) m_i m_j H(x_i - x_j)}.
\end{equation}
By $\Delta H=-J+\tilde J$,
\begin{eqnarray}
\mathbb{A}_0 X & = &     \big\vert \log \epsilon \big\vert^{-2}
  \sum_{i,j \in I_q}{ \Delta  H (x_i
- x_j) m_i m_j (d(m_i) + d(m_j))}  \\
   & = & -  \big\vert \log \epsilon \big\vert^{-2}  \sum_{i,j \in
   I_q}{ J (x_i - x_j) m_i m_j (d(m_i) + d(m_j))} \nonumber \\
    &  & +   \big\vert \log \epsilon \big\vert^{-2}  \sum_{i,j \in
   I_q}{ \tilde{J} (x_i - x_j) m_i m_j (d(m_i) + d(m_j))} \nonumber
\\
   &=&:\mathbb{A}_0^1 X+\mathbb{A}_0^2 X
\end{eqnarray}
 From the non-positivity of $\mathbb{A}_0^1 X$,
the non-positivity of $\mathbb{A}_C X$, apparent from (\ref{locc}),
and the non-negativity of $X$, follows
\begin{equation}\label{lokc}
- \mathbb{E}_N \int_0^T \mathbb{A}_0^1 X(t) dt
- \mathbb{E}_N \int_0^T \mathbb{A}_C X(t) dt  \leq \mathbb{E}_N X(0)+
  \mathbb{E}_N \int_0^T \mathbb{A}_0^2 X(t) dt .
\end{equation}
 $\Box$\\

\subsection{Bounds on functionals of $u_{n,m}$}\label{sttwo}
We will verify the assertions presented in (\ref{jayeye}). The
following lemma provides the bounds on the behaviour of the functions
$\big\{ u^\epsilon_{n,m}: \mathbb{R}^2 \to [0,\infty): (n,m) \in
\mathbb{N} \big\}$ and other functions that
will be used in this section.
  We choose the constant $R_0$ so that $V(x)=0$ whenever $|x|\ge R_0$.
Recall that $k=\sum_n nh_n$.
\begin{lemma}\label{phiest} There exists a collection of constants $C:
\mathbb{N}^2 \to (0,\infty)$ for which
the following bounds hold. \\
\begin{itemize}
\item for $x \in \mathbb{R}^2$ satisfying $|x|\le 2R_0\e$,
$\big\vert u_{n,m}^{\eps}(x) \big\vert \leq C_{n,m} |\log \e|$, and
for all $x \in \mathbb{R}^2$, 
$\big\vert u_{n,m}^{\eps}(x) \big\vert \leq C_{n,m} |\log
\vert x \vert |$.
\item for $x \in \mathbb{R}^2$,
$\big\vert \nabla u_{n,m}^{\eps} (x) \big\vert \leq
C_{n,m} \min \Big\{\frac{1}{\vert x \vert},\frac 1{\e}\Big\}$.
\item for $x \in \mathbb{R}^2$,
\begin{equation}\label{lept}
    \Big\vert u^\epsilon_{n,m} \big( x + z \big) - u^\epsilon_{n,m} \big( x
    \big) \Big\vert \leq C_{n,m}\vert z \vert \min\left\{
   \frac{1} {\vert x \vert},\frac 1\e\right\}
\end{equation}
\item for $x \in \mathbb{R}^2$ satisfying $\vert x \vert \geq   \max \big\{ 2
\vert z \vert + R_0 \epsilon , 2 R_0 \epsilon \big\}$,
\begin{equation}\label{lepta}
    \Big\vert \nabla u^\epsilon_{n,m} \big( x + z \big) - \nabla
u^\epsilon_{n,m} \big( x
    \big) \Big\vert \leq \frac{C_{n,m} \vert z \vert}{\vert x \vert^2}.
\end{equation}
\item let $H=H_{n,m}: \mathbb{R}^2 \times \R^2\to [0,\infty)$ be given by
\begin{displaymath}
    H (x;z) = \frac{-1}{2 \pi} \int_{\vert x - y \vert \leq
    1}{\log {\vert x - y \vert}} 
    u_{n,m}^\epsilon \big( y+z\big) 1 \! \! 1 \big\{ \vert y
    \vert \leq \rho \big\} dy .
\end{displaymath}
Then,
\begin{displaymath}
\sup_{|z|\le 1}\int{ H (x_1-x_2;z)k(x_1)k(x_2)dx_1
dx_2} \leq C_{n,m} \rho^2 |\log \rho|  .
\end{displaymath}
\item Let $\hat{H} \big( =
\hat{H}_{n,m} \big): \mathbb{R}^2\times \R^2 \to
[0,\infty)$ be given by
\begin{displaymath}
  \hat{H} (x;z)  = \frac{-1}{2\pi} \int_{\vert x - y \vert \leq
    1} {\log {\vert x - y \vert}} \  \Big\vert \nabla
    u_{n,m}^\epsilon \big( y+z\big) \Big\vert 1 \! \! 1 \big\{ \vert y
    \vert \leq \rho \big\} dy.
\end{displaymath}
Then, for every $z$ with $|z|\le 1$,
\begin{equation}\label{lepto}
\int{ \hat{H} (x_1-x_2;z)k(x_1)k(x_2)dx_1
dx_2} \leq C_{n,m} \big( \rho + \vert z \vert \big)  .
\end{equation}
\item for $(x,z) \in \mathbb{R}^2\times \R^2$, let $L (x;z)$ be given by
\begin{displaymath}
    \frac{-1}{2 \pi} \int_{\vert x - y \vert \leq 1} \log 
  {\vert x - y \vert} \left[\big|\log|1-|y+z||\big|+1\right]
   1 \! \! 1 \big\{ 1 - \rho \leq 
\vert y + z  \vert \leq 1 + \rho \big\} dy.
\end{displaymath}
Then we have the bound
\begin{equation}\label{eq2.20}
L(x;z)\le C\rho\big(\log\rho\big)^2 .
\end{equation}
\item 
for any positive integers $n$ and $m$ and a nonnegative smooth function
$\overline J$ of compact support, there exists
a constant $C_{n,m}(\overline J)$ such that, for any given $z \in
\mathbb{R}^2$, the function
$A_{n_1,n_2,n_3}^{\epsilon}: \mathbb{R}^6 \to [0,\infty)$ defined
by
\begin{eqnarray}
    & & \big( d(n_1 )\Delta_{x_1} + d(n_2) \Delta_{x_2} + d(n_3)
\Delta_{x_3} \big)
    A_{n_1,n_2,n_3}^{\epsilon} (x_1,x_2,x_3) \nonumber \\
    & & \qquad = -   u_{n,m}^\e(x_1 - x_3 + z)  V_{\eps}
    \big( x_1 - x_2 \big) 1\!\!1 \big\{ |x_1-x_3|\le \rho \big\}
\overline{J}(x_3) \nonumber
\end{eqnarray}
satisfies
\begin{eqnarray}
    & &\sum_{n_1,n_2,n_3}\int_{\mathbb{R}^6}{ A_{n_1,n_2,n_3}^{\epsilon} \big( 
x_1,x_2,x_3 \big)
    h_{n_1}(x_1)h_{n_2}(x_2)h_{n_3}(x_3) d x_1
    d x_2 d x_3  } \nonumber\\
    & & \qquad \qquad\leq C_{n,m}(\overline J) \parsize(\rho+|z|) \log
    \big( \rho + \vert z \vert \big) .
\label{leptu}
\end{eqnarray}
\end{itemize}
\end{lemma}
{\bf Proof}
Throughout the proof, we write $u^\e$ for the function $u_{n,m}^\e$ and
$\tau$ for the constant $\a(n,m)/(d(n)+d(m))$. The dependence of the
constants on $n$ and $m$ arises from that of $\tau$, and is also omitted.
The first part of the Lemma is a straightforward consequence
of our results in Section 3. As a consequence of Theorems 3.1-3.2 and Lemma 3.1 
we know that there exists a  constant $c_1\in (0,1)$ such that for
small $\e$ and $y$ satisfying
$|y|\le 2R_0$,
\[
\log \e \le u^\e(\e y)\le c_1\log\e,
\]
or equivalently,
\begin{equation}\label{eq2.30}
0\le u^\e(\e y)|\log \e|^{-1}+1\le 1-c_1.
\end{equation}
From (2.1),
we learn that for $x$ satisfying $|x|\ge 2\e R_0$,
$$
{\Lambda}_\e\log(|x|/2)\le u^\e(x)\le \log(2|x|){\Lambda}_\e,
$$
where
$$
{\Lambda}_\e=\frac 1{2\pi}\tau  \int \Big( V_{\epsilon} (y ) \ue_{n,m} ( y ) +  V^{\epsilon}( y ) \Big) dy.
$$
From this and (2.22) we learn that there are
two positive constants $k_1$ and $k_2$ such that if $|x|\ge 2\e R_0$,
then
$$
k_1\log(|x|/2))\le u^\e(x)\le k_2 \log(2|x|).
$$ 
To prove the second part of the lemma,
recall firstly that
\begin{displaymath}
    u^{\epsilon} \big( x \big) = \frac \tau {2\pi}\int_{\mathbb{R}^2}{
    \log{\vert x - y \vert}  \Big( u^{\epsilon}(y)
    V_{\epsilon}(y) + V^{\epsilon}(y) \Big) dy}.
\end{displaymath}
As a result,
\begin{equation}\label{wopdo}
    \nabla \ue(x) = \frac \tau {2\pi}\int_{\mathbb{R}^2}{
     \frac{x-y}{\big\vert x - y \big\vert^2}
     \Big( \ue(y)V_{\eps}(y) + V^{\eps}(y) \Big)  dy}.
\end{equation}
If $\vert x \vert\ge 2R_0 \e$, then $\big\vert x - y \big\vert \geq
\vert x \vert / 2$, and
\begin{displaymath}
    \bigg\vert \frac{x - y}{\vert x - y \vert^2} \bigg\vert =
    \frac{1}{\vert x - y \vert} \le \frac{2}{\vert x \vert},
\end{displaymath}
implying that 
\begin{displaymath}
    \big\vert \nabla \ue(x) \big\vert \leq \frac{ 2{\Lambda}_\e}{\vert x \vert}\le 
    \frac 3{|x|} ,
    \end{displaymath}
for small $\e$. If $\vert x \vert\le 2R_0\e$, then
we use (2.22) to deduce
\begin{displaymath}
    \big\vert \nabla \ue(x) \big\vert \le c_1
  \int_{|x-y|\le 3R_0\e}{\frac{1}{\vert x - y \vert}
  V^\e(y)dy}  \leq c_2 \eps^{-2}
  \int_0^{3R_0\eps}{dr} \leq c_3\e^{-1}.
\end{displaymath}
Thus,
\begin{displaymath}
\big\vert \nabla \ue (x) \big\vert \leq C \min \Big\{\frac 1\e, \frac{1}{\vert x
\vert} \Big\},
\end{displaymath}
as claimed in the second part of the lemma.

To prove the third part of the lemma, note that
\begin{eqnarray}
    & & \Big\vert u^\epsilon \big( x + z \big) - u^\epsilon \big( x
    \big) \Big\vert
    \nonumber \\
    & \leq & \frac 1{2\pi}  \int_{\mathbb{R}^2}
    \bigg\vert \log \frac{\vert x - y + z \vert}{\vert x - y
  \vert}   \bigg\vert
  \bigg( u^{\epsilon}(y) V_{\epsilon}(y) + V^{\epsilon}(y) \bigg) dy
  \nonumber \\
  & \leq & \frac 1{2\pi} \int_{\mathbb{R}^2} 
  \bigg\vert \log \frac{\vert x - y + z \vert}{\vert x - y
  \vert}  \bigg\vert V^{\epsilon}(y) dy,
  \label{dtc}
\end{eqnarray}
the latter inequality by means of (2.22).
From this and the elementary inequalities
\begin{equation}\label{rpp}
   \log \frac{\vert x - y + z \vert}{\vert x - y  \vert} -1
    \leq  \frac{\vert x-y+z
  \vert}{\vert x - y  \vert} \leq \frac{ \vert z \vert}{\vert x - y
  \vert},
\end{equation}
we deduce that
\[
\Big\vert u^\epsilon \big( x + z \big) - u^\epsilon \big( x
    \big) \Big\vert
     \leq  \frac 1{2\pi} |z| \int_{\mathbb{R}^2}
     \frac 1{2\pi} \int_{\mathbb{R}^2} 
  \frac{ V^{\epsilon}(y)}{|x-y|} dy
 \]
 We now use this and argue as in the proof of the second part of the lemma 
 to deduce the third part of the lemma.
 
In seeking to prove the fourth part of the lemma, note that
\begin{displaymath}
    \frac{x + z - y}{\big\vert x + z - y \big\vert^2} - \frac{x -
    y}{\big\vert x - y \big\vert^2} = \frac{\big\vert x - y \big\vert^2 \big(x + z - y \big)
     - \big\vert x + z - y
    \big\vert^2\big( x - y \big) }{ \big\vert x + z - y \big\vert^2 \big\vert x - y \big\vert^2 }.
\end{displaymath}
Note that, for any $a \in \mathbb{R}^2$,
\begin{eqnarray}
    \Big\vert \vert a \vert^2\big( a + z \big)  - \vert a + z \vert^2 a
    \Big\vert & \leq & \vert a \vert \Big\vert \vert a + z \vert^2 -
    \vert a \vert^2 \Big\vert + \vert z \vert \vert a \vert^2 \nonumber
    \\
    & \leq & c_1 \vert z \vert \vert a \vert^2, \label{cule}
\end{eqnarray}
so long as $|z|\le|a|$.
Note that since by our assumption $|x|\ge 2|z|+R_0\e$,
we have that $|x-y|\ge 2|z|$.
We may apply (\ref{cule}) with the choice $a = x - y$ to the formula
(\ref{wopdo}), hereby obtaining
\begin{displaymath}
     \Big\vert  \nabla\ue ( x + z ) - \nabla\ue (x)
     \Big\vert \leq C \vert z \vert \int_{\mathbb{R}^2}{
     \frac{V^{\e}(y)}{\big\vert x + z - y \big\vert^2}} dy,
\end{displaymath}
where we used (2.22).
 From the inequality $\vert x \vert \geq \max \big\{ 2 \vert z \vert +
R_0 \epsilon , 2 R_0 \epsilon \big\}$, we deduce that  $\big\vert x +
z - y \big\vert \geq \vert x -y \vert/2$ and $\vert x - y \vert \geq
\vert x \vert/2$. We conclude that
\begin{displaymath}
     \Big\vert  \nabla\ue ( x + z ) - \nabla\ue (x)
     \Big\vert \leq  \frac{C \vert z \vert}{\vert x \vert^2} ,
\end{displaymath}
as required.

To prove the fifth part of the lemma, note that
\begin{displaymath}
    |H (x;z)| \leq C \int_{|x-y|\le 1}{ |\log |x-y|| \log \vert y + z
    \vert\vert
1 \! \! 1 \big\{ \vert y \vert \leq \rho \big\} dy} \, ,
\end{displaymath}
by the first part of the lemma. Hence,
\begin{eqnarray}
    & & \int{|H(x_1-x_2;z)|k(x_1)k(x_2)dx_1dx_2} \nonumber \\
& \leq & C \int_{\vert y \vert \leq
    \rho}{  \big|\log \vert y + z
    \vert \big|   \int_{|x_1-x_2-y|\le 1}\big| \log
    {\vert x_1 - x_2 - y \vert} \big| k(x_1) k(x_2) dx_1
    dx_2 dy}
    \nonumber \\
    & \leq &  C \int_{\vert
    y \vert \leq \rho}{\big|\log \vert y \vert\big| dy} \nonumber \\
   & \leq & C \int_0^{\rho}{r |\log r| d r}
     \leq C \rho^2 |\log \rho|  \nonumber,
\end{eqnarray}
where in the second inequality, we used our first assumption on the initial data and
 the fact that if $\rho+|z|\le 1$, then the expression
 $\int_{\vert y \vert \leq \rho}|{\big|\log \vert
y+z \vert\big| dy}|$  is maximized as a function of $z \in \mathbb{R}^2$ when
$z = 0$.
This establishes the third part of the lemma.

To prove the sixth part of the lemma,
note that, by the first part,
\begin{displaymath}
    \hat{H} (x;z) \leq
  \frac{1}{2\pi} \int_{\vert x - y \vert \leq 1}{\big|
\log {\vert x - y \vert} \big| \frac{1}{\vert y + z
\vert}  1 \! \! 1 \big\{ \vert y \vert \leq \rho \big\} }dy.
\end{displaymath}
It follows that
\begin{eqnarray}
    & &  \int{|\hat H(x_1-x_2;z)|k(x_1)k(x_2)dx_1dx_2} \nonumber \\
& \leq & c_1 \int_{\vert y \vert \leq
    \rho}{ \frac{1}{\vert y + z  \vert} \int_{\vert x_1 - x_2 - y \vert
    \leq 1}{ \big|\log {\vert x_1 - x_2 - y \vert} \big|
    k(x_1)k(x_2)dx_1dx_2}dy} \nonumber \\
     & \le & c_2    \int_{\vert y \vert \leq
    \rho}{\frac{1}{\vert y + z  \vert} dy} 
    \le c_2    \int_{\vert y+z \vert \leq
    \rho+|z|}{\frac{1}{\vert y + z  \vert} dy}\nonumber \\
    & \leq & c_3  \big( \rho + \vert z \vert \big).
\end{eqnarray}
We have deduced (\ref{lepto}).

As for the seventh part of the lemma, first observe that $L(x;z)=L(x+z;0)$.
Hence we only need to verify (\ref{eq2.20}) when $z=0$. In this case we 
divide the domain of integration into the sets $|x-y|\le |1-|y||$ 
and $|1-|y||\le |x-y|\le 1$. Hence,
\begin{eqnarray}
  L(x;0) & \leq & \frac{1}{2\pi} \int_{\vert x - y \vert \leq \rho} 
\left[\big|\log {\vert x - y \vert} \big|+\big|\log {\vert x - y \vert} \big|^2
\right] dy \nonumber\\
& & +
 \frac{1}{2\pi} \int_{\big|1-|y|\big|\le \vert x - y \vert \le 1 } 
\left[\big|\log {|1-\vert y \vert}| \big|+\big|\log {|1-\vert y \vert}| \big|^2 
\right]1\!\!1\big(|y|\in (1-\rho,1+\rho)\big) dy 
\nonumber\\
 & \leq & \int_0^{\rho}{r \left[|\log r|+|\log r|^2\right] dr} \, + \,  
 \int_{1-\rho}^{1+\rho}\left[\big|\log|1-r|\big|+
 \big|\log|1-r|\big|^2\right]rdr\nonumber \\
  & \leq & C  \rho (\log \rho)^2 , \nonumber
\end{eqnarray}
establishing (\ref{eq2.20}).

As for the eighth part of the lemma, let us write $J(a)$ for
$1\!\!1 \big\{ |a|\le \rho \big\}$ and define the quantity $I$
according to
    \begin{eqnarray}
    I & = & c_0(6)\sum_{n_1,n_2,n_3} \g(n_1,n_2)\int_{} dx_1 dx_2 dx_3
    \int_{\mathbb{R}^{6}}  \Big( \frac{\vert x_1 - z' \vert^2}{d(n_1)}
+ \frac{\vert x_2 -
y \vert^2}{d(n_2)} + \frac{\vert x_3 - y' \vert^2}{d(n_3)}
\Big)^{-2}  \nonumber \\
    & & \qquad h_{n_1}(x_1)h_{n_2}(x_2)h_{n_3}(x_3) u^\e( z' - y' +
    z)  V_{\eps} ( z' - y ) J (z'-y')
\overline{J}(y') dz' dy dy'. \nonumber
\end{eqnarray}
We write
\[
I=c_0(6) \int_{\mathbb{R}^{6}}  u^\e(z' - y' +
    z)  V_{\eps} ( z' - y ) J (z'-y')
    \overline{J}(y') G(z',y,y')
    dz' dy dy',
    \]
    where $G(z',y,y')$ is given by
    \[
    \sum_{n_1,n_2,n_3 \g(n_1,n_2)\in \mathbb{N}}
    \int
    \bigg( \frac{\vert x_1 - z' \vert^2}{d(n_1)}  + \frac{\vert x_2 -
y \vert^2}{d(n_2)} + \frac{\vert x_3 - y' \vert^2}{d(n_3)}
\bigg)^{-2} h_{n_1}(x_1)h_{n_2}(x_2)h_{n_3}(x_3) dx_1 dx_2 dx_3.
\]
Using the elementary inequality $abc\le (a^2+b^2+c^2)^{3/2}$ we
deduce that $G(z',y,y')$ is at most
\begin{eqnarray}
    & &
\sum_{n_1,n_2,n_3 \in \mathbb{N}} \g(n_1,n_2){\left(d(n_1)d(n_2)d(n_3)\right)}^{2/3}
\int
    {\vert x_1 - z' \vert}^{-4/3}\ {\vert x_2 -
y \vert}^{-4/3} {\vert x_3 - y' \vert}^{-4/3} \nonumber \\
    & & \qquad \qquad \qquad h_{n_1}(x_1)h_{n_2}(x_2)h_{n_3}(x_3) dx_1
dx_2 dx_3.\nonumber
\end{eqnarray}
 From our assumptions on $h_n$ we deduce that $G\in L^{\i}_{loc}$. Hence,
\begin{displaymath}
    I  \leq   C \int_{\mathbb{R}^{6}}  u^\e(z' - y' +
    z) V_{\eps} ( z' - y ) J (z'-y')
    \overline{J}(y') dz' dy dy'.
\end{displaymath}
Note that, for fixed $z' \in \mathbb{R}^3$,
\begin{displaymath}
    \int_{\mathbb{R}^2}{V \Big( \frac{z' - y}{\epsilon}\Big) dy} =
\epsilon^2.
\end{displaymath}
Thus,
    \begin{eqnarray}
    I & \leq &  C \parsize \int_{\mathbb{R}^4}  u^\e(z' - y' +
    z)  J (z'-y') \overline{J}(y') dz' dy' \nonumber \\
    & \leq & C \parsize \int_{K} \int_{\mathbb{R}^2} u^\e(z' - y' +
    z) J(z'-y')dz' dy'  \nonumber \\
     & \leq & C \parsize \int_{K} dy' \int_{\mathbb{R}^2}\big| \log
     \vert z' - y' + z\vert\big| J(z'-y')dz' \\
     &\leq& C \parsize  \int_{|a|\le \rho+|z|}\big| \log \vert a \vert \big|da
  \leq C \parsize(\rho+|z|)^2 \big|\log \big( \rho + \vert z \vert \big)\big| , \nonumber
\end{eqnarray}
where $K \subseteq \mathbb{R}^2$ denotes a compact set containing the
support $\overline{J}$, and where we made use of the first
part of the lemma in the third inequality. This is the bound stated in
    (\ref{leptu}). $\Box$
\subsection{Estimating the terms}\label{stthree}
\subsubsection{The case of $H_{14}$ and $H_4$}
The estimate of
$\mathbb{E}_N \int_{0}^{T}{\vert H_4 (t) \vert dt}$ is derived in an
identical fashion to that of
$\mathbb{E}_N \int_{0}^{T}{\vert H_{14} (t) \vert dt}$.
Note that
\begin{eqnarray}
    \mathbb{E}_N \Big\vert \int_0^T H_{14} (t) dt \Big\vert & \leq &
    C \sparsize \mathbb{E}_N \int_0^T dt \sum_{i,j \in
    I_q}{ \Big\vert \ue \big( x_i - x_j + z \big) - \ue
    \big( x_i - x_j \big) \Big\vert } \nonumber \\
     & & \qquad \quad 1\!\!1 \Big\{ m_i = M_1, m_j = M_2 \Big\},
     \nonumber
\end{eqnarray}
where the constant $C$ depends on the $L^{\infty}$ bounds satisifed by
$J,\overline{J}$ and their time derivatives.
Hence
$$
  \mathbb{E}_N \Big\vert \int_0^T H_{14} (t) dt \Big\vert
  \leq K_1 + K_2,
$$
where $K_1$ is given by
\begin{eqnarray}
    & & C \sparsize \mathbb{E}_N \int_0^T{} dt\ \sum_{i,j \in I_q:
\vert x_i - x_j \vert >
\rho}{\Big\vert \ue(
x_i - x_j + z ) - \ue(
x_i - x_j  ) }\Big\vert  \nonumber \\
    & &  \qquad \qquad \qquad \qquad \qquad \qquad 1 \! \! 1 \{ m_i =
M_1 \} 1 \! \!  1 \{ m_j =
M_2 \} \nonumber
\end{eqnarray}
and $K_2$ is given by
\begin{eqnarray}
    & &  \sparsize \mathbb{E}_N \int_0^T{} dt\ \sum_{i,j \in I_q:
\vert x_i - x_j \vert \leq
\rho}{\Big\vert \ue(
x_i - x_j + z ) - \ue(
x_i - x_j  ) } \Big\vert \nonumber \\
    & &   \qquad \qquad \qquad \qquad  \qquad \qquad 1 \! \! 1 \{ m_i =
M_1 \} 1 \! \!  1 \{ m_j = M_2 \}  . \nonumber
\end{eqnarray}
Firstly, we treat $K_1$. Note that
\[
K_1   \leq   \frac{C \vert z \vert \sparsize}{\rho}
\mathbb{E}_N \int_0^T{}
dt\ \sum_{i,j \in I_q: \vert x_i - x_j \vert > \rho}{ 1 \! \! 1 \{ m_i = M_1
\}   1 \! \! 1 \{ m_j = M_2 \} }  \leq   \frac{C
    \vert z \vert Z^2 }{\rho}  , 
\]
where the first inequality follows from the third part of Lemma
\ref{phiest}, and the final one from the initial number of particles
$N$ equals $Z|\log\e|$.

We now treat the term $K_2$. By writing,
\begin{eqnarray}
    K_2 & \leq & C \sparsize \mathbb{E}_N \int_0^T{} dt\ \sum_{i,j \in
I_q: \vert x_i - x_j \vert \leq
\rho}{
    \bigg[ \big|\ue \big(
x_i - x_j + z \big)\big| + \big|\ue \big(
x_i - x_j  \big)\big| \bigg] }  \nonumber \\
    & &   \qquad \qquad \qquad \qquad  \qquad \qquad 1 \! \! 1 \{ m_i =
M_1 \} 1 \! \!  1 \{ m_j = M_2 \}  , \nonumber
\end{eqnarray}
we obtain an expression on the right-hand-side which may be bounded by
applying Lemma \ref{lemtwop}. As a result we can write $K_2\le K_{21}+K_{22}$
where $K_{21}$ and $K_{22}$ represent the first and the second term on the right-hand-side in
Lemma \ref{lemtwop}. For $K_{21}$, the relevant estimate is provided by
the fifth part of Lemma \ref{phiest}, with a bound of
$C\rho^2 |\log \rho|$.
To bound the  term $K_{22}$, note that, with the function $J$ in Lemma
\ref{lemtwop} chosen to be $J(x) = u^{\epsilon}(x + z)1\!\!1(|x|\le \rho)$, we have that
\begin{eqnarray}
  |\tilde{J}(x) |& = &\left| \int_{\vert y \vert = 1}{u^{\epsilon}\big( x + y +
  z \big) 1 \! \! 1 \big\{ \vert x + y \vert \leq \rho \big\} }S(dy)\right| \nonumber
  \\
  & \leq & C \int_{\vert y \vert = 1}{\big|\log \vert x + y + z \vert\big| 1 \!
  \! 1 \big\{ \vert x + y \vert \leq \rho \big\}}S(dy) \nonumber 
\end{eqnarray}
 Let us assume that  $|z|\le \rho$ ans set $a=x+z$. We  then have
\begin{eqnarray*}
  |\tilde{J}(x) |
  & \leq & C \int_{\vert y \vert = 1}{|\log \vert a +y  \vert| 1 \!
  \! 1 \big\{ \vert a + y \vert \leq 2\rho \big\}}S(dy)  \\
  & \leq & 2C 
  \int_0^{c_2\rho}{|\log (c_1 r)| dr} \leq c_3 \rho |\log \rho| ,
\end{eqnarray*}
 where  for the second inequality we used that fact that the conditions
$|a+y|\le \rho, \ |y|=1$ mean that the point $y$ belongs to an arc on the unit circle 
of center $\sigma=a/|a|$ and length $2c_2 \rho$, and if the lenght of the arc $\sigma y$
is $r$, then $|a+y|\ge c_1 r$ for positive constants $c_1 $ and $c_2$. 

In this way, we find that $\tilde J$ is uniformly bounded by
$C\rho |\log\rho|$ and this in turn implies that the  term $K_{22}$ is bounded above by
\begin{displaymath}
  C \big\vert \log \eps \big\vert^{-2} \rho |\log\rho|
  \mathbb{E}_N \int_0^T dt \Big\vert \Big\{ (i,j ) \in I_q^2 : m_i =
  M_1, m_j = M_2 \Big\} \Big\vert \leq C Z^2 \rho |\log\rho|,
\end{displaymath}
the latter inequality following from $N \leq Z \big\vert \log \eps
\big\vert$. We find that
\begin{displaymath}
  K_2 \leq C  \rho |\log\rho| \, .
\end{displaymath}
Hence,
\begin{displaymath}
\mathbb{E}_N \left|\int_0^T{}  H_{13}(t) \right| dt  \leq  K_1 + K_2
\leq    \frac{C
    \vert z \vert  }{\rho}   +  C \rho |\log\rho|  .
\end{displaymath}
Setting $\rho = \vert z \vert^{1/2}$, we find that
\begin{displaymath}
\mathbb{E}_N \left|\int_0^T{}  H_{13}(t) \right| dt  \leq C  \vert z
\vert^{\frac{1}{2}} \big|\log {\vert z \vert} \big| .
\end{displaymath}
\subsubsection{The cases of $H_2$ and $H_3$}
The estimate of
$\mathbb{E}_N \int_{0}^{T}{\vert H_3 (t) \vert dt}$ is derived in an
identical fashion to that of
$\mathbb{E}_N \int_{0}^{T}{\vert H_2 (t) \vert dt}$.
Picking $\rho \in \mathbb{R}$ that satisfies $\rho \geq \max \big\{ 2 \vert z
\vert + R_0 \epsilon , 2 R_0 \epsilon \big\}$,
we write
\begin{equation}\label{eq2.28}
\int_{0}^{T}{\mathbb{E}_N{\vert H_2 (t) \vert} dt} \leq R_1 + R_2,
\end{equation}
where
\begin{eqnarray}
R_1 & = &  \sparsize \mathbb{E}_N{} \int_{0}^{T}{\sum_{i,j \in
I_q: \vert x_i
- x_j \vert > \rho}{d(m_i)  \left|\overline{J}(x_j,m_j,t)\right|}} \nonumber \\
    & & \qquad \quad \Big\vert \ue_x ( x_i - x_j + z ) -
    \ue_x ( x_i - x_j ) \Big\vert
      \big\vert J_x \big( x_i,m_i,t \big) \big\vert dt
, \nonumber
\end{eqnarray}
and
\begin{eqnarray}
R_2 & = &  \sparsize \mathbb{E}_N \int_{0}^{T}{\sum_{i,j \in
I_q: \vert x_i
- x_j \vert \leq \rho}{d(m_i)  \big|\overline{J}(x_j,m_j,t)\big|}} \nonumber \\
& & \qquad \quad
    \Big\vert \ue_x( x_i - x_j + z ) -
    \ue_x ( x_i - x_j ) \Big\vert
     \big\vert J_x \big( x_i,m_i,t \big) \big\vert dt . \nonumber
\end{eqnarray}
Firstly, we examine the sum $R_1$.
Recalling that we consider test functions $J$ and $\overline{J}$
respectively supported on particles of mass $M_1$ and $M_2$,
\begin{displaymath}
R_1  \leq  C \sparsize   \frac{\vert
z \vert}{\rho^2} d(M_1)
\Big\vert \{ (i,j) \in I_q^2: m_i = M_1, m_j = M_2 \} \Big\vert,
\end{displaymath}
where the lower bound on $\rho$ allowed us to apply the fourth part of
Lemma \ref{phiest}. Thus, 
\begin{equation}\label{eq2.29}
R_1 \leq C \frac {\vert z \vert}{\rho^2}.
\end{equation}

Secondly, we bound the sum $R_2$. Note that
\begin{eqnarray}
R_2 & \leq &  \sparsize  \vert \vert J_x \vert \vert \ \vert \vert
\overline{J} \vert \vert \ d(M_1)
    \  \mathbb{E}_N \int_0^T{}dt \sum_{i,j \in I_q: \vert x_i - x_j
\vert \leq \rho} \nonumber \\
    & & \qquad   \Big\vert
    \ue_x(x_i - x_j + z ) - \ue_x(
x_i - x_j  ) \Big\vert  1 \! \! 1 \{ m_i = M_1, m_j = M_2 \} \nonumber \\
& \leq &  C \sparsize \ \mathbb{E}_N \int_0^T{}dt \sum_{i,j \in I_q}{m_i
m_j \big( d(m_i) + d(m_j) \big)} \nonumber \\
    & & \qquad \qquad \quad  \bigg[ \Big\vert
    \ue_x(x_i - x_j + z ) \Big\vert + \Big\vert \ue_x(
x_i - x_j  ) \Big\vert   \bigg] 1 \! \! 1 \Big\{ \big\vert x_i - x_j
\big\vert \leq \rho \Big\}  ,  \label{fhg}
\end{eqnarray}
where $\|\cdot \|$ denotes the $L^\i$ norm and
the constant $C$ depends on the test functions $J$ and $\overline{J}$.
The expression (\ref{fhg}) is written in a form to which Lemma \ref{lemtwop}
may be applied. Doing so yields
\begin{eqnarray}
R_2 & \leq & C \sparsize   \mathbb{E}_N  \sum_{i,j
\in I_q}{m_i m_j \Big( \hat{H} (x_i - x_j;0) + \hat{H} (x_i - x_j;z)
\Big) } \nonumber \\
   & & \, + \,  C \sparsize \ \mathbb{E}_N \int_0^T{}dt \sum_{i,j \in I_q}{m_i
m_j \big( d(m_i) + d(m_j) \big)}  \nonumber\\
    & & \qquad \qquad \quad  \Big[
    \tilde{J}\big( x_i - x_j ;z \big)  +
    \tilde{J}\big( x_i - x_j ;0 \big)   \Big]    \nonumber\\
    &=& : R_{21}+R_{22}\label{eq2.31}
\end{eqnarray}
where the function $ \hat{H} $ appears in the sixth part of Lemma
\ref{phiest}, and where $\tilde{J}: \mathbb{R}^2\times \R^2 \to [0,\infty)$ in
this case is given by
\begin{displaymath}
  \tilde{J}(x;z) = \frac 1{2\pi}\int_{\vert y \vert = 1}{\big\vert u_x^{\eps}(x+y+z)
  \big\vert 1 \! \! 1 \big\{ \vert x + y \vert \leq \rho \big\} dS(y) }.
\end{displaymath}

 From the sixth part of Lemma \ref{phiest} and our assumptions on the initial
data, we deduce that 
\begin{equation}\label{eq2.32}
R_{21}\le C(\rho + \vert z \vert).
\end{equation}

It follows from the second part of Lemma \ref{phiest} that the 
function $\tilde{J}$ satisfies the bound 
$$
|\tilde{J}(x;z)|\le C\int_{\vert y \vert = 1}\frac 1{\big\vert x+y+z \big\vert}
 1 \! \! 1 \big\{ \vert x + y \vert \leq \rho \big\} dS(y)  .
 $$
By our assumption, we certainly have  $|z|\le \rho/2$. Hence,
$$
|\tilde{J}(x;z)|\le C\int_{\vert y \vert = 1}\frac 1{\big\vert a+y
  \big\vert }1 \! \! 1 \big\{ \vert a + y \vert \leq 2\rho \big\} dS(y)  ,
$$
for $a=x+z$. Note that there exists a positive constant $c_1$
such that the conditions $|a+y|\le 2\rho,\ |y|=1$ mean
that $y\in \G $, where $\G$ is an arc of the unit circle with the center 
$\sigma=-a/|a|$. 
It is not hard to show that there exist positive constants $c_1$ and $c_2$ such
that \[
c_1\big(\ell+\big|1-|a|\big|\big)\le |a+y|\le c_2\big(\ell+\big|1-|a|\big|\big),
\]
where $\ell$ denotes the length
 of the arc from
$\sigma=-a/|a|$ to $y$ on the unit circle. From this we deduce 
\begin{eqnarray*}
|\tilde{J}(x;z)|&\le &1 \! \! 1 \big\{ 1-\rho\le 
 \vert a  \vert \leq 1+\rho \big\}\int
 1\!\!1\left({c_1c_2^{-1}\big|1-|a|\big|}\le \ell\le {2c_1^{-1}\rho}\right)
 \frac {c_3}\ell d\ell\\ 
  &\le& c_3 1 \! \! 1 \big\{ 1-\rho\le  \vert x+z  \vert \leq 1+\rho \big\}
  \left(\big|\log|1-|x+z||\big|+\log\big(c_2c_1^{-2}\big)\right) 
\end{eqnarray*}

Using this bound on $\tilde J$, and then applying Lemma \ref{lemtwop}, we learn that 
the term $R_{22}$ is bounded above by $R_{221}+R_{222}$, where
\begin{displaymath}
  R_{221} = C \big\vert \log \epsilon \big\vert^{-2} \mathbb{E}_N 
\sum_{i,j \in I_{q(0)}}{\Big( L \big( x_i - x_j;0 \big) + L \big( 
x_i - x_j;z \big) \Big) m_i m_j},
\end{displaymath}
and,
\begin{displaymath}
  R_{222}= 2C \big\vert \log \epsilon \big\vert^{-2} \mathbb{E}_N
    \int_0^T \Big( \tilde{L} \big( x_i - x_j \big) + \tilde{L}
\big( x_i - x_j+z \big) \Big) m_i m_j \big( d(m_i) + d(m_j) \big) dt,
\end{displaymath}
where the function $L: \mathbb{R}^2 \to [0,\infty)$ appears in the 
seventh part of Lemma \ref{phiest},  and where the function 
$\tilde{L} : \mathbb{R}^2 \to [0,\infty)$ is given by
\begin{displaymath}
  \tilde{L}(x) =\frac 1{2\pi} \int_{\vert y \vert = 1}{ 1 \! \! 1 \big\{ 1 - \rho 
\leq\vert x  + y \vert \leq 1 + \rho \big\}  \left[\big|\log|1-|x+y||\big|+1
\right] d S(y) }
\end{displaymath}
By the seventh part of Lemma \ref{phiest} and our assumption
on the initial total density $k=\sum_n nh_n\in L^1$ we obtain 
\begin{equation}\label{eq2.33}
R_{221} \leq C\rho 
(\log \rho )^2.
\end{equation}

We decompose $R_{222}=R_{2221}+R_{2222}$ where 
\begin{displaymath}
  R_{222r}= 2 \big\vert \log \epsilon \big\vert^{-2} \mathbb{E}_N
    \int_0^T \Big( {L}_r \big( x_i - x_j \big) + {L}_r
    \big( x_i - x_j+z \big) \Big) m_i m_j \big( d(m_i) + d(m_j) \big) dt,
\end{displaymath}
for $r=1,2$, where $L_1(x)=\tilde L(x)1\!\!1(|x|\ge \rho^{1/4})$ and
$L_2(x)=\tilde L(x)1\!\!1(|x|\le\rho^{1/4})$.
Let us now analyse the behaviour of the function $\tilde L$. 
First observe that since the Lebesgue measure on the circle is rotationally invariant, 
we have that the function $\tilde L$ is radially symmetric. Because of this,
let us assume that $x=(a,0)$ and $y=(cos\theta,\sin\theta)$. 
Note that
$$
|x+y|^2=1+a^2+2a\cos\theta.
$$
As a result, the condition $|x+y|\in (1-2\rho,1+2\rho)$ for $\rho\le 1$ implies
that $|a^2+2a\cos\theta|\le 8\rho$. 
Let us first examine the case $|a|=|x|\ge  \rho^{1/4}$. In this case we have that 
$|a+2\cos\theta|\le 8\rho^{3/4}$. This condition is not satisfied
unless $|a|\le 2$. In that case, choose $\theta_0\in [0,\pi]$ such that $a+2\cos\theta_0
=0$. Note that for small$\rho$, the set $\{\theta:|a+2\cos\theta|\le 8\rho^{3/4}\}$
is a union of two disjoint $\theta$--intervals
about the points $\theta_0$ and $\theta_1=2\pi-\theta_0$.
 We now argue that there exists a positive 
constant $c_1$ such that the length of these intervals is bounded above by $c_1\rho^{1/2}$.
To see this, first observe that the condition $|a|\ge \rho^{1/4}$ implies that for a
positive constant $c_2$ we have that
$|\theta_0|,|\theta_0-\pi|\ge 2c_2\rho^{1/4}$. This and $|a+2\cos \theta|\le 8 \rho^{3/4}$
implies that we also have $|\theta|,|\theta-\pi|\ge c_2\rho^{1/4}$
 provided that $\rho$ is sufficiently small. On the other hand since
 \[
 2\cos \theta+a=-2(\sin \tau) (\theta-\theta_0),
 \]
 for some $\tau$ between $\theta$ and $\theta_0$, we deduce that for some positive constant
 $c_1$, we have that $|\theta-\theta_0|\le c_1\rho^{1/2}$. 
 Also, there exist positive constants $c_3$ and $c_4$ such that if $\theta $ is 
 close to $\theta_r$, then
$$
\big|1-|x+y|\big|\ge c_3|a^2+2a\cos\theta|\ge c_4 \sqrt \rho |\theta-\theta_r|,
$$
for $r=0$ or $1$.
From this we learn that if $|x|\ge  \rho^{1/4}$, then
the term $|\tilde L(x)|$ is bounded above by
$$
 \frac 1{2\pi}\int_{\theta_0-c_1\sqrt\rho}^{\theta_0+c_1\sqrt\rho}
\left[|\log \left(c_4\sqrt \rho |\theta-\theta_0|\right)|+1\right]d\theta
+\frac 1{2\pi}\int_{\theta_1-c_1\sqrt\rho}^{\theta_1+c_1\sqrt\rho}
\left[|\log \left(c_4\sqrt \rho |\theta-\theta_1|\right)|+1\right]d\theta.
$$
As a result, $|\tilde L(x)|\le C\sqrt\rho|\log\rho|$. This in turn implies that
\begin{equation}
R_{2221}\le C\sqrt \rho|\log\rho| 
\end{equation}
We now turn to $R_{2222}$. For this, observe that the support of the 
function $L_2$ is contained in the set of points $x$ for which $|x|\le 
 \rho^{1/4}$. Note that if $|a|\le  \rho^{1/4}$, then
 we can find a positive constant $c_5$ such that
  $|\theta_0-\pi/2|\le c_5\rho^{1/4}$,
$|\theta_1-3\pi/2|\le c_5\rho^{1/4}$, where
$\theta_1=2\pi-\theta_0$. Forthermore, we can find a positive constant $c_6$ such that 
if $\theta\in (0,\pi)$, then 
\begin{equation}\label{eq2.37}
\big|a^2+2a\cos\theta\big|=\big|2a(\cos\theta_0-\cos\theta)\big|
=2\left|a\sin\frac {\theta+\theta_0}2\sin\frac{\theta-\theta_0}2\right|
\ge c_6 |a(\theta-\theta_0)|.
\end{equation}
The same is true if $\theta\in [\pi,2\pi]$ (use $\theta_1$ in place of $\theta_0$
in (\ref{eq2.37}).)
 As a result,
$$
\big|1-|x+y|\big|\ge c_3 |a^2+2a\cos\theta|\ge c_3c_6 |a(\theta-\theta_0)|.
$$
From this we learn that indeed
$$
L_2(x)\le C\big|\log |x|\big|1\!\!1\big\{|x|\le \rho^{1/4}\big\}.
$$
To bound $R_{2222}$, let us apply Lemma 2.2 one more time to write $R_{2222}\le
R_{22221}+R_{22222}$, where 
\begin{displaymath}
  R_{22221} = C \big\vert \log \epsilon \big\vert^{-2} \mathbb{E}_N 
\sum_{i,j \in I_{q(0)}}{\Big( \Gamma \big( x_i - x_j \big) + \Gamma \big( 
x_i - x_j+z \big) \Big) m_i m_j}
\end{displaymath}
and
\begin{displaymath}
  R_{22222}=  \big\vert \log \epsilon \big\vert^{-2} \mathbb{E}_N
    \int_0^T \Big( \tilde{\Gamma} \big( x_i - x_j \big) + \tilde{\Gamma}
\big( x_i - x_j+z \big) \Big) m_i m_j \big( d(m_i) + d(m_j) \big) dt,
\end{displaymath}
where the function $\Gamma: \mathbb{R}^2 \to [0,\infty)$ is very similar to 
the function $H$ that appeared in the 
fifth part of Lemma \ref{phiest} (except that $\rho$ in the definition
$H$ is replaced with $\rho^{1/4}$),
 and where the function 
$\tilde{\Gamma} : \mathbb{R}^2 \to [0,\infty)$ is given by
\begin{displaymath}
  \tilde{\Gamma}(x) = \int_{\vert y \vert = 1}{ 1 \! \! 1 \big\{ \vert
   x  + y \vert \leq \rho^{1/4} \big\}  \big|\log|x+y|\big| d S(y) }.
\end{displaymath}
As in the fifth part of Lemma 2.4 we show
\begin{equation}\label{eq2.35}
R_{22221}\le C\rho|\log\rho|.
\end{equation}
In just the same way that we bounded
$\tilde J$ in the subsection 3.3.1, we can readily show that $\tilde \Gamma(x)\le
C\sqrt\rho|\log \rho|$. This in turn implies
\begin{equation}\label{eq2.36}
R_{22222}\le C\rho^{1/4}|\log \rho|.
\end{equation}
Putting all the pieces together we learn from 
(\ref{eq2.28})--(\ref{eq2.29}) and (\ref{eq2.31})--(\ref{eq2.36}) that 
\begin{displaymath}
\int_0^T{\mathbb{E}_N{\vert H_2(t) \vert} dt}  
    \leq C \left[ \frac{\vert z \vert}{\rho^2}  + \rho + \vert z \vert
    +{\rho}^{1/4}
     \big|\log\rho \big|,\right]
\end{displaymath}
for $\rho\le 1$. By making the choice $\rho =  \vert z
\vert^{\frac{4}{9}} $, we find
that
\begin{displaymath}
\int_0^T{\mathbb{E}_N{\vert H_2(t) \vert} dt}  \leq  C  \vert z
\vert^{\frac{1}{9}}\big|\log|z|\big|.
\end{displaymath}
\subsubsection{The case of $G_z(1) - G_0(1)$}
We now estimate the term
$$
\int_{0}^{T}{ \mathbb{E}_N \big\vert G_z(1) - G_0(1) \big\vert (t) dt}.
$$
To ease the notation, we do not display the dependence of $J$ and 
$\overline J$ on the variable $t$.
Note that
\begin{equation}\label{eq2.40}
\int_{0}^{T}{ \mathbb{E}_N \big\vert G_z(1) - G_0(1) \big\vert (t) dt}
\leq \sum_{i=1}^{8}{D_i},
\end{equation}
where
\begin{eqnarray}
D_1 & = & \frac 12\mathbb{E}_N \int_0^T{} dt \sum_{k,l \in
I_q}{\micp(m_k,m_l)V_{\epsilon}(x_k - x_l)} \big\vert
J(x_k,m_k)\big\vert   \nonumber\\
& & \quad \epsilon^{2(d-2)} \sum_{i \in I_q}{
    \big\vert\overline{J}(x_i,m_i)\big\vert   \Big\vert
\ue(x_k - x_i + z) - \ue(x_k - x_i )
\Big\vert} , \nonumber
\end{eqnarray}
each of the other seven terms on the right-hand-side of (\ref{eq2.40})
differing from $D_1$ only in an inessential way.
Given this. the estimates involved for each of the eight cases are in
essence identical, and we examine only the case of $D_1$. We write $D_1 =
D^1 + D^2$,
where we have decomposed the inner $i$-indexed sum according to the
respective index sets
$$ \{ i \in I_q, i \not= k,l, \vert x_k - x_i \vert > \rho \} \
\textrm{and} \  \{ i \in I_q, i \not= k,l, \vert x_k - x_i \vert \leq \rho \}.
$$
By the second part of Lemma \ref{phiest}, we have that
$$
D^1  \leq  \frac{C \vert z \vert
\parsize}{\rho} \mathbb{E}_N \int_0^T{} dt \sum_{k,l \in
I_{q}}{\alpha(m_k,m_l) V_{\epsilon}(x_k - x_l)},  \nonumber
$$
where we have also used the fact that the test functions $J$ and
$\overline{J}$ are each supported on the set of particles of
respective masses $M_1$ and $M_2$, and the fact that the total number
of particles living at any given time is bounded
above by $Z \big\vert \log \epsilon \big\vert$. From the bound on the 
collision that is
provided by Lemma
\ref{lembc}, it follows that
$$ D^1 \leq  \frac{C \vert z \vert}{\rho} .$$

Note that $D^2$ is bounded above by
\begin{eqnarray}
    &  & C \mathbb{E}_N \int_0^T{} dt \sum_{k,l \in
I_{q}}{\alpha(m_k,m_l) V_{\epsilon}(x_k - x_l) } 1 \! \! 1 \{ m_k =
M_1 \} \sparsize \sum_{i \in I_q}\nonumber \\
& & \quad {} 1 \! \! 1 \{ \vert x_i - x_k
\vert \leq \rho\} 1 \! \! 1 \{ m_i = M_2 \} \Big\vert
\ue(x_k - x_i + z) - \ue(x_k - x_i ) \Big\vert \big|J(x_i,m_i)\big| 
\nonumber \\
    & \leq & C \mathbb{E}_N \int_0^T{} dt \sum_{k,l \in
I_{q}}{\alpha(m_k,m_l)V_{\epsilon}(x_k - x_l)} 1 \! \! 1 \{ m_k = M_1
\} \sparsize \sum_{i \in I_q}\nonumber \\
& & \quad {} 1 \! \! 1 \{ \vert x_i - x_k
\vert \leq \rho \} 1 \! \! 1 \{ m_i = M_2 \}  \Big[
|\ue(x_k - x_i + z)| + |\ue(x_k - x_i )| \Big]\big|J(x_i,m_i)\big| . \nonumber
\end{eqnarray}
Note that the last expectation is bounded by
Lemma \ref{lemthreep} because, by our assumption on $\a$, we can find
$\g$ such that
$\a\le \g$, with  $\g$ satisfying the assumption of Lemma 2.3. The upper
bound provided by this Lemma in this
particular application is computed in the last part of Lemma
\ref{phiest}. We find that $D^2 \leq C (\rho+|z|) \log \big( \rho +
\vert z \vert \big)$.

Combining these estimates yields
$$
D_3 \leq D^1 + D^2 \leq C \frac{\vert z \vert}{\rho} +
  C  (\rho+|z|) \big|\log \big( \rho + \vert z \vert \big) \big|.
$$
Making the choice $\rho = \vert z \vert^{\frac{1}{2}} $ leads to the
inequality $D_3 \leq   \vert z \vert^{\frac{1}{2}}\big|\log \vert z \vert\big|$.
Since each of the cases of $\big\{ D_i : i \in \{ 1,\ldots, 8 \} \big\}$ may
be treated by a nearly verbatim
proof, we deduce that
$$
\int_{0}^{T}{ \mathbb{E}_N \big\vert G_z(1) - G_0(1) \big\vert (t) dt}
\leq   C \vert z \vert^{\frac{1}{2}} \log \vert z \vert .
$$
\subsubsection{The case of $G_z(2)$} Recall that
\begin{displaymath}
G_z(2)  =  - \sparsize \sum_{k,l \in
I_{q}}{ \alpha(m_k,m_l) V_{\epsilon}(x_k - x_l)} \ue(x_k -
x_l + z ) J(x_k,m_k)
    \overline{J}(x_l,m_l).
\end{displaymath}
If $k,l \in I_q$ satisfy $V_\epsilon \big( x_k - x_l \big) \not= 0$,
then $\vert x_k -
x_l \vert \leq R_0 \epsilon$, and so
$$
\vert x_k - x_l + z \vert \geq \vert z \vert - R_0 \epsilon
\geq \vert z \vert/2  ,
$$
provided that $|z|\ge 2R_0\e$. This implies that
$$
    \big|\ue(x_k - x_l + z )\big| \leq
    C \Big\vert \log \vert x_k - x_l + z \vert \Big\vert
  \leq C \Big\vert \log \vert z \vert \Big\vert ,
$$
where in the first inequality, we used the first part of Lemma
\ref{phiest} (restated).
Applying this bound, and using the fact that the test functions $J$
and $\overline{J}$ have compact support, we find that
$$
\int_0^T{} \mathbb{E}_N \vert G_z (2) \vert dt \leq C \Big\vert \log
\vert z \vert \Big\vert
\sparsize  \mathbb{E}_N
\int_0^T{} \sum_{k,l \in I_q}{ \alpha(m_k,m_l) V_{\epsilon}(x_k - x_l)} dt
$$
whose right-hand-side is bounded above by $C \big\vert \log \vert z
\vert \big\vert \parsize$,
according to Lemma \ref{lembc}. That is,
$$
\int_0^T{\vert G_z(2) \vert dt} \leq C  \big\vert \log \vert z
\vert \big\vert  \parsize.
$$
\subsubsection{The case of $\mathbb{E}_N \vert X_z - X_0 \vert $}\label{lak}

We now turn to $\mathbb{E}_N \vert X_z - X_0 \vert $. Assume that 
$|z|\ge R_0\e$.
Using the  second part of Lemma 2.4, we have that
\[
\mathbb{E}_N \vert X_z - X_0 (0)\vert \le  C|z|
\iint_{L^2} h_{M_1}(x)h_{M_2}(y)|x-y|^{-1}
dxdy,
\]
where $L$ is a bounded set that contains the support of $J$ and $\overline J$.
Using our second assumption on the initial data  $h_n$ we obtain the 
bound $C|z|$ for  $\mathbb{E}_N \vert X_z - X_0 \vert(0)$.

\begin{subsection}{The martingale term}\label{stfour}
This section is devoted to proving the estimate (\ref{martest}).
Note that
$$
M_z(T) = X_z \big( q(T),T \big) - X_z \big( q(0),0 \big) -
\int_{0}^{T}{\left(\frac \partial {\partial t}+\mathbb{L} \right)X_z 
\big( q(t),t) \big) dt}
$$
is a martingale which satisfies
$$
\mathbb{E}_N \Big[  M_z(T)^2 \Big] = \sum_{i=1}^{3}{}
\mathbb{E}_N\int_{0}^{T}{A_i \big( q(t)
,t
\big) dt},
$$
where $A_1(q,t)$ and $A_2(q,t)$ are respectively set equal to
\begin{eqnarray}
&  & 2  \big\vert \log \eps \big\vert^{-4} \sum_{i \in I_q, m_i  = M_1}{d(M_1)}
\nonumber \\
    & & \qquad \qquad \quad
    \bigg| \nabla_{x_i}
\sum_{j \in I_q, m_j = M_2}{\ue(x_i - x_j + z)
J(x_i,M_1,t) \overline{J}(x_j,M_2,t) } \bigg|^2, \nonumber
\end{eqnarray}
and
\begin{eqnarray}
&  & 2  \big\vert \log \eps \big\vert^{-4}
   \sum_{i \in I_q, m_i  = M_2}{d(M_2)}
\nonumber \\
    & & \qquad \qquad \quad
    \bigg| \nabla_{x_i}
\sum_{j \in I_q, m_j = M_1}{\ue(x_j - x_i + z)
    J(x_j,M_1,t) \overline{J}(x_i,M_2,t) } \bigg|^2, \nonumber
\end{eqnarray}
while $A_3$ is given by
\begin{eqnarray}
&  &  \frac 12 \big\vert \log \eps \big\vert^{-4}  \sum_{i,j \in I_q}{\alpha(m_i,m_j)
  V_{\eps} ( x_i - x_j ) } \label{aatwo} \\
     & & \qquad \quad \bigg\{ \sum_{k \in I_q}{} \Big[
     \frac{m_i}{m_i + m_j} \ue (x_i - x_k+z )  J(x_i,m_i +
     m_j) \overline{J}(x_k,m_k) \nonumber \\
    & & \qquad \qquad \qquad + \
     \frac{m_i}{m_i + m_j} \ue (x_k - x_i + z)  J(x_k,m_k)
\overline{J}(x_i,m_i + m_j) \nonumber \\
    & & \qquad \qquad \qquad + \
     \frac{m_j}{m_i + m_j} \ue (x_j - x_k + z) J(x_j,m_i +
     m_j) \overline{J}(x_k,m_k) \nonumber \\
    & & \qquad \qquad \qquad + \
     \frac{m_j}{m_i + m_j} \ue (x_k - x_j + z)  J(x_k,m_k)
\overline{J}(x_j,m_i + m_j) \nonumber \\
    & & \qquad \qquad \qquad - \  \ue (x_i - x_k + z)
     J(x_i,m_i) \overline{J}(x_k,m_k) 
    - \  \ue (x_k - x_i + z)
     J(x_k,m_k) \overline{J}(x_i,m_i) \nonumber \\
      & & \qquad \qquad \qquad - \   \ue (x_j
     - x_k + z)   J(x_j,m_j) \overline{J}(x_k,m_k) 
   - \   \ue (x_k
     - x_j + z)  J(x_k,m_k) \overline{J}(x_j,m_j) \Big] \nonumber \\
     & & \qquad \qquad \qquad \qquad - \  \ue (x_i
     - x_j + z) J(x_i,m_i) \overline{J}(x_j,m_j)
      \bigg\}^2. \nonumber
\end{eqnarray}
Recall that, by our convention, we do not display the dependence of 
$J$ and $\overline J$ on the
$t$-variable. To bound these terms, we require two variants of Lemma 2.3 :
\begin{lemma}\label{tpln}
There exists a collection of constants $C: \mathbb{N}^2 \to
(0,\infty)$ such that,
for any continuous functions $t,v,a_1,a_2,a_3: \mathbb{R}^2 \to 
[0,\infty)$ and any $z\in \R^2$,
\begin{eqnarray}
& & \mathbb{E}_N\int_0^T{}dt \sum_{i,j,k \in I_q(t)}{}\g(m_i,m_j)t ( x_i -
x_j+z ) v ( x_i - x_k+z )
\nonumber \\
    & & \qquad \qquad \qquad a_1(x_i)
a_2(x_j) a_3(x_k) 1 \! \! 1 \big\{ m_i = n_1 ,
m_k = n_3 \big\}   \nonumber \\
& & \qquad \qquad \quad  \leq C_{n_1,n_3} |\log\e|^{3}\sum_{n_2}
\mathbb{E}_N
\sum_{i,j,k \in I_{q(0)} }A_{n_1,n_2,n_3}^{\epsilon}(x_i,x_j,x_k)  , \nonumber
\end{eqnarray}
where $A_{n_1,n_2,n_3}^{\epsilon}: \mathbb{R}^6 \to [0,\infty)$ is
given by
\begin{eqnarray}
     A_{n_1,n_2,n_3}^{\epsilon}(x_1,x_2,x_3) & = &
c_0(6)\g(n_1,n_2)\int_{\mathbb{R}^6} \bigg(
\frac{\vert x_1 - z' \vert^2}{d(n_1)} + \frac{\vert x_2 - y
\vert^2}{d(n_2)} + \frac{\vert x_3 - y' \vert^2}{d(n_3)}
\bigg)^{-2} \nonumber \\
    & & \quad  t (  z' - y+z ) v
( z' - y'+z ) a_1(z') a_2(y) a_3(y') dz' dy dy, \nonumber
\end{eqnarray}
with $\g$ as in Lemma 2.3.
\end{lemma}
\begin{lemma}\label{fpln}
There exists a collection of constants $C: \mathbb{N}^3 \to
[0,\infty)$ such that,
for any $z\in \R^d$,
  any continuous functions $v,w: \mathbb{R}^2 \to [0,\infty)$ and
another $(a_1,a_2,a_3): \mathbb{R}^6 \to [0,\infty)$,
\begin{eqnarray}
& & \mathbb{E}_N\int_0^T{}dt \sum_{k,l,i,j \in I_q}{}\g(n_i,n_j)V_{\e}( {x_i -
x_j}) v ( x_i - x_k+z ) w
( x_i - x_l+z ) a_1( x_i)a_2(x_k)a_3(x_l) \nonumber \\
    & & \qquad \qquad \qquad  1 \! \! 1 \big\{ m_i = n_1 ,  m_k = n_3,
m_l = n_3 \big\}   \nonumber \\
& & \qquad \quad \leq C_{n_1,n_3,n_3} \big\vert \log \epsilon \big\vert^4
\sum_{n_2}\mathbb{E}_N \sum_{i,j,k,l
\in I_{q(0)}}{B_{m_i,m_j,m_k,m_l}^{\epsilon}(x_i,x_j,x_k,x_l)} \nonumber
\\
& & \qquad \qquad \qquad \qquad 1 \! \! 1 \big\{ m_i \leq n_1 , m_j
\leq n_2, m_k \leq n_3,
m_l \leq n_3 \big\} , \nonumber
\end{eqnarray}
where  $B_{n_1,n_2,n_3,n_3}^{\epsilon}: \mathbb{R}^8 \to [0,\infty)$ is
given by
\begin{eqnarray}
    & & B_{n_1,n_2,n_3,n_3}^{\epsilon} \big( x_1,x_2,x_3,x_3 \big)
    \nonumber \\
    & = & c_0(8)\int_{\mathbb{R}^{8}} \bigg(
\frac{\vert x_1 - \hat{z} \vert^2}{d(n_1)} + \frac{\vert x_2 - z'
\vert^2}{d(n_2)} + \frac{\vert x_3 - y \vert^2}{d(n_3)} +
\frac{\vert x_3 - y' \vert^2}{d(n_3)}
\bigg)^{-3} \nonumber \\
    & & \qquad \g(n_1,n_2)V_\e( {\hat{z} - z'}) v ( \hat{z} - y+z )
  w ( z' - y'+z ) a_1( \hat{z})a_2(y)a_3( y') d \hat{z}
dz' dy dy', \nonumber
\end{eqnarray}
with the function $\g:\mathbb{N}^2\to (0,\i)$ satisfying
\begin{displaymath}
n_2 \g \big( n_1, n_2 + n_3 \big) \max{ \Big\{ 1, \Big[ \frac{d (
n_2 + n_3 )}{d(n_2)}
\Big]^3 \Big\}} \leq \big( n_2 + n_3 \big) \g(n_1,n_2).
\end{displaymath}
\end{lemma}

The proof of Lemma 2.5 is identical to that of Lemma 2.3. The proof
of Lemma 2.6 is very similar
to the proof of Lemma 2.3 and is omitted.

We now bound the three terms. Of the first two, we treat only $A_1$,
the other being bounded by an identical argument.
By multiplying out the brackets appearing in the definition of $A_1$,
we obtain that this quantity is bounded above by
\begin{eqnarray}
    & & C \big\vert \log \epsilon \big\vert^{-4}
   \sum_{i,j,k \in I_q} \big\vert  \nabla \ue
    \big\vert ( x_i - x_j + z )   \big\vert
    \nabla \ue \big\vert (  x_i - x_k + z )  J^2(x_i,m_i) \nonumber \\
& &  \qquad \qquad
    \big|\overline{J}(x_j,m_j)\big| \big|\overline{J}(x_k,m_k)\big|  1 
\! \! 1 \big\{
m_i = M_1 , m_j = m_k = M_2
    \big\}\\
    & + & C  \big\vert \log \epsilon \big\vert^{-4}
  \sum_{i,j,k \in I_q}  \left|\ue ( x_i - x_j + z  )\right|\left|   \ue ( x_i -
x_k + z  )\right| |\nabla J(x_i,m_i)|^2 \nonumber \\
    & &  \qquad \qquad
    \big|\overline{J}(x_j,m_j)\big| \big|\overline{J}(x_k,m_k)\big|  1 
\! \! 1 \big\{
m_i = M_1 , m_j = m_k = M_2
    \big\}  . \label{yhj}
\end{eqnarray}
Let us assume that $z=0$ because this will not affect our arguments.
We are required to bound the quantity appearing in
the statement of Lemma \ref{tpln}, for each of the following cases:
\begin{equation}\label{twocases}
    \big( t,v,a_1,a_2,a_3 \big)  \in  \Big\{ \big( \vert\nabla \ue
\vert , \vert \nabla \ue \vert,
    J^2, \big|\overline{J}\big|,  \big|\overline{J}\big|\big) \ , \   \big(
  \ue,  \ue,  \vert \nabla J \vert^2 ,  \big|\overline{J}\big|, 
\big|\overline{J}
\big|
    \big) \Big\}.
\end{equation}
Recall that each  of the test functions $J$, $\overline{J}$, and
their gradients,
is assumed to be uniformly bounded with compact support. To each of
the two cases, Lemma \ref{tpln} applies. For either of them, the
right-hand-side of the inequality in Lemma \ref{tpln} may be written as a
finite sum of the expectations appearing there, with the sum being
taken over triples of given masses $n_1,n_2$ and $n_3$. Such an expectation is
bounded above by
\begin{eqnarray}
    & & C  \big\vert \log \epsilon \big\vert^{-4}  \sum_{n_2 \in \mathbb{N}}
    \int \int_{K^3} \bigg(
\frac{\vert x_1 - z' \vert^2}{d(n_1)} + \frac{\vert x_2 - y
\vert^2}{d(n_2)} + \frac{\vert x_3 - y' \vert^2}{d(n_3)}
\bigg)^{-2} \label{yun} \\
    & & \qquad \qquad t ( z' - y )
v ( z' - y' ) h_{n_1}(x_1)h_{n_2}(x_2)h_{n_3}(x_3) dz'
dy dy'  d x_1 d x_2 d x_3 , \nonumber
\end{eqnarray}
where $K=\{x:|x|\le \ell\} \subseteq \mathbb{R}^2$ is chosen to
contain the support of $J$ and $\overline{J}$. As in Section 3.4 of \cite{HR}, we can use our
bounds in the first two parts of Lemma 2.4 and repeat the proof of the eighth part
of Lemma 2.4 to obtain
\begin{equation}\label{eq2.44}
\mathbb{E}_N\int_0^T\left[A_1 \big( q(t) \big)+A_2 \big( q(t) \big)
\right] dt \leq C \parsize .
\end{equation}
We must treat the third term, $A_3$.
An application of the inequality
$$
\Big( a_1 + \ldots + a_n \Big)^2 \leq n \Big( a_1^2 + \ldots + a_n^2 \Big)
$$
to the bound on $A_3$ provided in (\ref{aatwo}) implies that
\begin{equation}\label{aathr}
    A_3(q) \leq \frac 92   \big\vert \log \epsilon \big\vert^{-4}
    \sum_{i,j \in I_q}{\alpha (m_i,m_j)
    V_{\eps} ( x_i - x_j ) \bigg[
    \sum_{n=1}^{8}{\Big( \sum_{k \in I_q}Y_n \Big)^2 } + Y_9^2  \bigg]},
\end{equation}
where $Y_1$ is given by
\begin{displaymath}
     \frac{m_i}{m_i + m_j} \ue (x_i - x_k + z ) J(x_i,m_i +
     m_j) \overline{J}(x_k,m_k),
\end{displaymath}
and where $\big\{ Y_i : i \in \{ 2, \ldots, 8 \} \big\}$ denote the
other seven expressions in (\ref{aatwo}) that appear in a sum over $k
\in I_q$, while $Y_9$ denotes the last term in (\ref{aatwo})
that does not appear in this sum.
There are nine cases to consider. The first eight are practically
identical, and we treat only the fifth. Note that
\begin{eqnarray}
    & &   \big\vert \log \epsilon \big\vert^{-4}
          \sum_{i,j \in I_q}{\alpha (m_i,m_j)  V_\e (
          x_i - x_j ) \Big(
    \sum_{k \in I_q}Y_5 \Big)^2 } \nonumber\\
    & =& C \big\vert \log \epsilon \big\vert^{-4} \sum_{i,j \in I_q}{\alpha(m_i,m_j)
    V_\e( x_i - x_j )} \nonumber \\
    & &  \qquad \Big[ \sum_{k,l \in I_q}{\ue ( x_i - x_k + z )
    \ue ( x_i - x_l + z )
    J^2(x_i,m_i) \overline{J}(x_k,m_k) \overline{J}(x_l,m_l) } \Big]. 
\label{zvc}
\end{eqnarray}
In the sum with indices involving $k,l \in I_q$, we permit the
possibility that these two may be equal, though they must be distinct
from each of $i$ and $j$ (which of course must themselves be distinct
by the overall convention).

Note that the expression (\ref{zvc}) appears in the statement of
Lemma \ref{fpln},
    provided that the choice
$$
\Big( v,w ,a_1,a_2,a_3\Big) = \Big(   \big|\ue\big|, \big|\ue\big|,
 J^2 ,\big| \overline{J}\big|,
    \big|\overline{J}\big|\Big)
$$
is made. Again we  set $z=0$ because this does not
affect the estimates.
Given that the support of each of the functions $a_1,a_2,a_3: 
\mathbb{R}^{6} \to
    [0,\infty)$ is bounded,
we must bound
\begin{eqnarray}
    & & \sum_{n_2 \in \mathbb{N}}\int  \int_{L^4} \Big(
\frac{\vert x_1 - \hat{z} \vert^2}{d(n_1)} + \frac{\vert x_2 - z'
\vert^2}{d(n_2)} + \frac{\vert x_3 - y \vert^2}{d(n_3)} +
\frac{\vert x_3 - y' \vert^2}{d(n_3)}
\Big)^{-3} V \Big( \frac{\hat{z} - z'}{\epsilon} \Big)\nonumber \\
    & & \qquad   \big|\ue ( \hat{z} - y )\big|
 \big| \ue ( z' - y' )\big| h_{n_1}(x_1)h_{n_2}(x_2)h_{n_3}(x_3)h_{n_3}(x_3)
    d \hat{z} dz' dy dy' d x_1 d x_2 d x_3  dx_3 , \nonumber
\end{eqnarray}
for  a compact set $L$. This expression is bounded above by
\begin{eqnarray}
     & & \int_{L^3}  V \Big( \frac{\hat{z} - z'}{\epsilon} \Big) \ue
  ( \hat{z} - y )  \ue ( z' - y' )  d \hat{z} dz' dy dy' \nonumber \\
    & & \qquad  \sum_{n_2 \in \mathbb{N}}{ \int_{K^3}  \bigg(
\frac{\vert x_1 - \hat{z} \vert^2}{d(n_1)} + \frac{\vert x_2 - z'
\vert^2}{d(n_2)} + \frac{\vert x_3 - y \vert^2}{d(n_3)} +
\frac{\vert x_3 - y' \vert^2}{d(n_3)}
    \bigg)^{-3}  }  \nonumber \\
    & & \qquad \qquad
    h_{n_1}(x_1)h_{n_2}(x_2)h_{n_3}(x_3)h_{n_3}(x_3)dx_1 dx_2 dx_3 
dx_3  , \nonumber
\end{eqnarray}
which is less than
$$
C  \int_{L^3}  V \Big( \frac{\hat{z} - z'}{\epsilon} \Big) \ue
  ( \hat{z} - y )  \ue ( z' - y' )  d \hat{z} dz' dy dy'.
$$
The proof of this follows the proof of the eighth part of Lemma 2.4;
we use the elementary
inequality $abcd\le (a^2+b^2+c^2+d^2)^2$ and the fact that the function
$$
\hat k(x)=\sum_nd(n)^{3/4}\int h_n(y)|x-y|^{-3/2}dy
$$
is locally bounded. Noting that the bound $\vert \ue(x) \vert \leq
\big\vert \log\vert x \vert \big\vert$ implies that
$$
\int_{L}{u^\e(\hat z-y )dy}
$$
is bounded above by a constant, we find that
$$
     \int_{L^3}  V \Big( \frac{\hat{z} - z'}{\epsilon} \Big)
\ue ( \hat{z} - y )  \ue ( z' - y' )  d\hat{z} dz' dy dy'
     \leq  C  \int_{L^2}  V \Big( \frac{\hat{z} -
z'}{\epsilon} \Big)   d \hat{z} dz'.
    $$
This is at most $ C\e^2$.  Applying Lemma \ref{fpln},
    we find that the contribution to
$$
\mathbb{E}_N{\int_{0}^{T}{A_3 \big(q(t)\big) dt}}
$$
arising from the fifth term in (\ref{aathr}) is at most
$$
C      \eps^{-2} 
\big\vert^{-1}    \eps^2
  = C \big\vert \log \eps \big\vert^{-1}.
$$

We now treat the ninth term, as they are classified in
(\ref{aathr}). It takes the form
\begin{eqnarray}
    & &    \big\vert \log \epsilon \big\vert^{-4}
      \sum_{i,j \in I_q}{\alpha(m_i,m_j)
     V_\e\Big( {x_i - x_j} \Big)} \nonumber \\
    & & \qquad \qquad \qquad
    \ue \big(x_i - x_j + z \big)^2
      J \big( x_i,m_i \big)^2 \overline{J} \big( x_j,m_j \big)^2. \nonumber
\end{eqnarray}
This is bounded above by
$$
      C \big\vert \log \epsilon \big\vert^{-3}
   \sum_{i,j \in I_q}{\alpha(m_i,m_j)
     V_\e\Big( {x_i - x_j} \Big)},
     $$
  because $u^\e\le C   \big\vert \log \epsilon \big\vert $ in the support
  of $V_\e$ by the first part of
  Lemma 2.4. The expected value of the integral on the interval of time $[0,T]$
of this last expression is bounded above by
$$
C   \big\vert \log \epsilon \big\vert^{-3}
  \mathbb{E}_N \int_0^T{dt} \sum_{i,j \in I_q}{\alpha(m_i,m_j)
    V_{\e}(x_i - x_j)}\le C   \big\vert \log \epsilon \big\vert^{-1}  ,
$$
where we used Lemma \ref{lembc} for the last inequality.
This completes the proof of (\ref{martest}).
\end{subsection}
\subsection{Using the estimates}\label{stfive}
The inequalities (\ref{abc}), (\ref{jayeye}) and (\ref{martest}) imply
that, for large $T$,
$$
\lim_{\vert z \vert \to 0}\limsup_{\epsilon \downarrow
0}\mathbb{E}_N   \bigg\vert\int_0^T H_{11} \big( t \big) dt  +
\int_0^T H_{13} \big( t \big) dt\bigg\vert dt = 0.
$$
That is,
\begin{eqnarray}
& &  \lim_{\vert z \vert \to 0}{\limsup_{\epsilon \downarrow
0}}    \sparsize \mathbb{E}_N  \bigg\vert
     \int_{0}^{T}{dt \sum_{i,j \in
    I_{q(t)}}{\alpha(m_i,m_j) J(x_i,m_i) \overline{J}(x_j,m_j)} } \nonumber \\
& & \qquad  \qquad  \Big[ V^{\epsilon} (x_i - x_j + z) -
V^{\epsilon} ( x_i - x_j )   + 
V_{\epsilon}(x_i - x_j + z) \ue(x_i - x_j + z)
    \Big] \ \bigg\vert     =\ 0. \nonumber
\end{eqnarray}
(Recall that we simply write $J(x_i,m_i)$ and $\overline J(x_i,m_i)$
for $J(x_i,m_i,t)$ and $\overline J(x_i,m_i,t)$.) This
implies that
\begin{eqnarray}
& &
    \sparsize \int_{0}^{T}  {\sum_{i,j \in
I_{q(t)}}{\alpha(m_i,m_j)
    V^{\epsilon} ( x_i - x_j ) J(x_i,m_i)}
    \overline{J}(x_j,m_j) } dt\nonumber \\
& = & \sparsize \int_{0}^{T}  {\sum_{i,j \in
I_{q(t)}}{\alpha(m_i,m_j)
    V^{\epsilon}( x_i - x_j + z ) J(x_i,m_i)
    \overline{J}(x_j,m_j)}}\nonumber \\
& & \qquad \qquad \qquad  \ \Big[ 1 +
     \parsize \ue(x_i - x_j
+ z)   \Big] dt  \ + \ Err_1(\epsilon,z)  , \label{otwdfw}
\end{eqnarray}
where $Err_1$ satisfies
\begin{equation}\label{lati}
    \lim_{\vert z \vert \to 0}{\limsup_{\epsilon \downarrow
    0}}{ \, \mathbb{E}_N \big\vert Err_1
    \big( \epsilon,z \big) \big\vert} = 0.
\end{equation}
By Theorem 3.2, the expression $1 + |\log\e|^{-1} u_{m_1,m_2}^\e \big( a
\big)$ is uniformly close to 
$$
 \left(1-\frac {{\tau}(m_1,m_2)}{2\pi+{\tau}(m_1,m_2)}\right),
$$
for $a$ satisfying $V^\e(a)\neq 0$ and ${\tau}(m_1,m_2)=\a(m_1,m_2)/(d(m_1)+
d(m_2))$. 
 Recalling from (\ref{qudef}) that
\begin{displaymath}
Q = \parsize \sum_{(i,j) \in
I_{q}}{\alpha(m_i,m_j) V_{\epsilon}(x_i - x_j) J(x_i,m_i)
    \overline{J}(x_j,m_j)}
\end{displaymath}
and writing
\begin{equation}\label{quzdef}
    \overline{Q} (z)  =  \parsize
    \sum_{i,j \in I_q}{\b(m_i,m_j) V^\e(
    {x_i - x_j + z
}) J(x_i,m_i)
    \overline{J}(x_j,m_j) },
\end{equation}
it follows from (\ref{otwdfw}) and Theorem 3.2 that
\begin{equation}\label{weakstoss}
\int_0^T  Q (t) dt  =
\int_0^T \overline{Q} (z)(t) dt  + Err_2\big(\epsilon,z \big),
\end{equation}
where $Err_2$ satisfies
\begin{displaymath}
    \lim_{\vert z \vert \to 0}{\limsup_{\epsilon \downarrow
    0}}{ \, \mathbb{E}_N \big\vert Err_2
    \big( \epsilon,z \big) \big\vert} = 0.
\end{displaymath}
From this, it is not hard to deduce that
\begin{eqnarray}
  \overline{Q}(z_2 - z_1) & = & |\log\e|^{-2}
    \sum_{i,j \in I_q}{\b(m_i,m_j) V^\e({x_i - x_j + z_2 - z_1 })}
     \label{strc} \\
   & & \qquad \qquad \qquad J(x_i - z_1,m_i)
    \overline{J}(x_j - z_2,m_j)  \, + \, Err(\epsilon,z_1,z_2), \nonumber
\end{eqnarray}
where
$$
\mathbb{E}_N \big\vert Err(\epsilon,z_1,z_2)\big\vert \leq
  C  \big( \vert z_1 \vert + \vert
    z_2 \vert \big).
$$
(See Section 3.5 of \cite{HR}.) By (\ref{weakstoss}) and (\ref{strc}),
\begin{eqnarray}
    & & \int_0^T  Q (t) dt  \nonumber \\
    & = & |\log\e|^{-2}\int_{0}^{T}dt
\sum_{i,j \in I_{q(t)}}\b(m_i,m_j) J(x_i - z_1,m_i)
\overline{J}(x_j - z_2,m_j)\nonumber \\
    & & \quad \int_{\mathbb{R}^2}{\int_{\mathbb{R}^2}{V^\e
    ( (x_i - z_1) - (x_j - z_2)) {\delta}^{-2}
\eta \Big( \frac{z_1}{\delta} \Big) {\delta}^{-2} \eta \Big(
\frac{z_2}{\delta} \Big) dz_1
dz_2 }} \, + \, Err_3\big(\epsilon,\delta\big) \nonumber \\
    &  =  &  |\log\e|^{-2} \int_{0}^{T} dt
\int_{\mathbb{R}^{2d}} d \omega_1 d \omega_2 \sum_{i,j
\in I_{q(t)}}{}V^\e( {\omega_1 -\omega_2})  \b(m_i,m_j)  J(\omega_1,m_i)
\overline{J}(\omega_2,m_j) \nonumber \\
    & & \qquad \quad
\delta^{-2} \eta \Big( \frac{x_i - \omega_1}{\delta} \Big)
\delta^{-2} \eta \Big( \frac{x_j - \omega_1}{\delta} \Big)
     \, + \, Err_3 \big( \epsilon, \delta \big) \nonumber \\
& = & \int_0^T{} dt \int_{\mathbb{R}^2} \int_{\mathbb{R}^2}
d\omega_1 d\omega_2 v^\e ( {\omega_1 -
\omega_2}) \b(M_1,M_2)  J(\omega_1,M_1)
\overline{J}(\omega_2,M_2) \nonumber \\
& & \qquad \quad  \bigg[ \parsize \sum_{i \in I_q: m_i = M_1}{
\delta^{-2} \eta \Big( \frac{x_i - \omega_1}{\delta} \Big) }
\bigg] \nonumber \\
& & \qquad \quad \bigg[ \parsize \sum_{j \in I_q; m_j = M_2}{
\delta^{-2} \eta \Big(\frac{x_j - \omega_2}{\delta} \Big)
} \bigg]  \, + \, Err_3 \big( \epsilon, \delta
\big) ,\nonumber
\end{eqnarray}
where $Err_3$ satisfies
$$
  \lim_{\delta \downarrow 0}\limsup_{\epsilon_\downarrow 0} \mathbb{E}_N
  \big\vert Err_3 \big( \epsilon, \delta \big) \big\vert = 0,
$$
and where in the last equality, we made use of the fact that the test
functions $J$ and $\overline{J}$ take non-zero values only on
particles of a given mass, respectively $M_1$ and $M_2$.
Thus,
\begin{eqnarray}
    & & \int_0^T  Q (t) dt  \nonumber \\
    & = &
\int_0^T{} dt \int_{\mathbb{R}^2}
d\omega  \b(M_1,M_2)  J(\omega,M_1)
\overline{J}(\omega,M_2)  \bigg[ \parsize \sum_{i \in I_q: m_i = M_1}{
\delta^{-2} \eta \Big( \frac{x_i - \omega}{\delta} \Big) }
\bigg] \nonumber \\
& & \qquad \quad \bigg[ \parsize \sum_{j \in I_q; m_j = M_2}{
\delta^{-2} \eta \Big(\frac{x_j - \omega}{\delta} \Big)
} \bigg]  \, + Err \big( \epsilon, \delta
\big) \nonumber
\end{eqnarray}
where $Err=Err_5+err$ with the function $err=O\big(\e \d^{-5}\big)$ also satisfies
$$
  \lim_{\delta \downarrow 0}\limsup_{\epsilon_\downarrow 0} \mathbb{E}_N
  \big\vert err \big( \epsilon, \delta \big) \big\vert = 0.
$$
This completes the proof of Proposition~\ref{szprop}.
\end{section}

\begin{section}{Potential theory}
\label{sec3}

The purpose of this section is twofold.  Firstly, we show the existence of the
function $u^\e$ that satisfies (\ref{eq3.1}). Secondly we evaluate
the limit of $u^\e |\log\e|^{-1}$ in the support of $V^\e$,
as $\e\to 0$. This limit was used in the evaluation of $\b$ in Section 2.5.
  We start with the
statements of the main
results of this section.  Let $V: {\R}^d \to {\R}$ be a continuous
function of compact support with $V \ge
0$ and $\ds{\int_{\mathbb{R}^2} V(x)dx = 1}$.  We also write $K_0$ 
for the topological closure
of $U_0$ where
\setcounter{equation}{0}
\begin{equation}
U_0 = \{x: V(x) \ne 0\}.
\end{equation}
Given a measure $\mu$, let us define
$$
\cG \mu(x)=\int \log |x-y| \mu(dy).
$$
When the measure $\mu$ is absolutely continuous with respect to the 
Lebesgue measure with a
density $g$, we simply write $\cG g$ for $\cG \mu$.

\begin{theorem}
\label{th3.1}
There exists a number ${\gamma}_0>0$ such that
for every ${\gamma} \in (0, {\gamma}_0)$ and  $a\in \R$, there exists a unique 
function $u \in
C^2({\R}^2)$ such that
$u(x)= O\big(\big|\log |x|\big|\big)$ as $|x|\to\i$ and
\begin{equation}
\label{eq3.2}
  u = {\gamma} \cG\big((u+a )V\big).
\end{equation}
Moreover $Z={\gamma} \int (u+a)Vdx\ne 0$ and $(u+a)Z^{-1}\ge 0$.
\end{theorem}

Recall that we are searching for a function $u^\e$ such that
$$
\Delta u^\e={\tau}(n,m)\left[V_\e u^\e+V^\e\right],
$$
where $\tau={\tau}(n,m)=\a(n,m)/(d(n)+d(m))$. For this it suffices to have
\begin{equation}\label{eq3.3}
u^\e=\cG \mu^\e,
\end{equation}
where $\mu^\e(dx)=\frac 1{2\pi} {\tau} (V_\e u^\e
+V^\e) dx$.
This can be rewritten as
$$
  u^\e=\frac 1{2\pi} {\tau}^\e\cG \left(V^\e\left[  u^\e+|\log\e|\right]\right),
$$
where ${\tau}^\e={\tau} |\log\e|^{-1}$. Evidently we can apply Theorem 3.1 to 
deduce the existence of the
function $u^\e$ for sufficiently small $\e$.

Our next theorem was used in the previous section for the evaluation of 
$\b$.

\begin{theorem}
\label{th3.2}
For every positive $k$,
\begin{displaymath}
\lim_{\e\to 0}\sup_{|x|\le k}\left|u^{\e}(\e x)|\log \e|^{-1}+\frac {\tau}{2\pi+{\tau}}\right|=0.
\end{displaymath}
\end{theorem}

\bigskip
\noindent
{\bf Proof of Theorem \ref{th3.1}} \\
\noindent{\bf Step 1.}  Let $J$ be a
bounded continuous
function with $J > 0$ and
$$
{\int_{|x| \ge 1}
J(x)\left(\log |x|\right)^{2}dx <
\i}.
$$
   Define
\begin{equation}
\label{eq3.5}
\cH = \left\{ u : \mbox{ $u$ is measurable and }\int_{\mathbb{R}^2} 
u^2(x) J(x) dx <
\i\right\}.
\end{equation}
We then define $\cF: \cH \to \cH$ by
$\cF(u)=\cG (uV)$.
Observe that $\cH$ is a Hilbert space with respect to
the inner product
\[
\<u,v\> = \int_{\mathbb{R}^2} u(x)v(x) J(x)dx.
\]
Let us verify that $\cF$ is a bounded operator.  To see this, write
\begin{equation}
\label{eq3.6}
\G(x) =  \int_{\mathbb{R}^2} \big |\log|x-y|\big |V(y)dy.
\end{equation}
When $|x|$ is sufficiently large, we have that $\G(x)\le 
 \log (2|x|)$ because
$0<\log|x-y|\le \log(2|x|)$ whenever $V(y)\ne 0$. Otherwise we have
$$
\G(x)\le c_0\int_{|x-y|\le c_1
}\big |\log|x-y|\big |dy\le c_2 ,
$$
for constants $c_0,c_1$ and $c_2$.
As a result
\begin{equation}\label{eq3.7}
\G(x)\le c+\log^+|x|
\end{equation}
for a constant $c$.  Also, we may use H\"older's inequality to assert
\begin{eqnarray}
\label{eq6.10}
(\cF(u)(x))^2 &\le &\left[ 
 \G(x) \int_{\mathbb{R}^2} \big |\log|x-y|\big |V(y)|u(y)| \frac
{dy}{\G(x)} \right]^2
\\
&\le & \G(x)\int_{\mathbb{R}^2} \big |\log|x-y|\big |V(y)u^2(y)dy. \nonumber
\end{eqnarray}
    From this  we deduce
\[
\int_{\mathbb{R}^2} (\cF(u)(x))^2J(x)dx \le \int_{\mathbb{R}^2} 
V(y)u^2(y) \left[
\int_{\mathbb{R}^2}
 \G(x)\big |\log|x-y|\big |J(x)dx\right]dy.
\]
If $V(y) \ne 0$ then $|y| \le R_0$ for a suitable $R_0$. Define
\[
I(y) = \int_{\mathbb{R}^2} \G(x)\big |\log|x-y|\big |J(x)dx,
\]
Note
\begin{eqnarray*}
I(y) &= &\int_{|x| \le 2R_0} + \int_{|x| > 2R_0}
\G(x)\big|\log|x-y|\big| J(x)dx \\
&\le &c_1 \int_{|x-y| \le 3R_0} \big |\log|x-y|\big |dx + c_1 \int_{|x| > 2R_0}
\big|\log|x|\big|^2 J(x)dx ,
\end{eqnarray*}
where for the second line we have used \eqref{eq3.7}.
From this and  our assumption on $J$ we deduce that
$\ds{\sup_{|y| \le R_0} I(y) < \i}$. As a result,
\[
\int_{\mathbb{R}^2}(\cF(u)(x))^2 J(x)dx \le c_1 \int_{\mathbb{R}^2} 
V(y)u^2(y)dy \le c_2  \int_{\mathbb{R}^2} 
u^2(y)J(y)dy
\]
because $V$ is of compact support and $J>0$.  This shows the 
boundedness of the operator $\cF:\cH\to \cH$.\\
\noindent{\bf Step 2.}  
Since the operator $\cF$ is bounded, the equation
\begin{displaymath}
(id - {\gamma}\cF)(u) = g
\end{displaymath}
has a solution, where $g(x) = -{\gamma}a\G(x)$ with $\G$ as in \eqref{eq3.6} 
and $id$ denotes
the identity transformation.
 Note that our assumption on $\G$
implies that $\G \in
\cH$ because of \eqref{eq3.7}.  
 So far we have shown the existence of a unique
solution $u \in \cH$ of
$u - {\tau} \cF(u) = g$.  From this and the H\"older continuity of $V$
we can readily show that in fact $u\in C^2$ and that $u$ is a
classical solution of
\begin{equation}
\Delta u=2\pi \gamma (u+a)V.
\end{equation}
(See for example Section 4.2 of \cite{GT}.)\\
\noindent{\bf Step 3.}  In this step we verify $Z\ne 0$.
Observe that $u=\cG \mu$ for a measure $\mu$ with a bounded support. From this
we can readily deduce
\begin{equation}
\label{eq3.9}
  u(x)=\mu\left(\R^2\right)\log|x|+O\left(|x|^{-1}\right),
\end{equation}
\begin{equation}
\label{eq3.10}
\nabla u(x)=
\mu\left(\R^2\right)\frac {x}{|x|^2}+O\left(|x|^{-2}\right).
\end{equation}
  We now
choose $R > R_0$ and use $\D u =2\pi{\gamma}  (u+a)V$ to write
\[
    \int_{|x| \le R} (u+a)\D u dx =2\pi {\gamma}\int_{|x| \le R} V(u+a)^2dx.
\]
After an integration by parts we obtain
\[
- \int_{|x| \le R} |\nabla u|^2dx +  \int_{|x| = R} (u+a) \nabla u
\cdot n dS =2\pi{\gamma}
\int_{|x| \le R} V(u+a)^2dx,
\]
where $\ds{n = \frac {x}{|x|}}$ is the normal vector and $dS$ is the
Lebesgue measure on $|x| = R$.  Now if $Z=\mu\big(\R^2\big)=0$,
then we can use \eqref{eq3.9}--\eqref{eq3.10} to deduce that
\[
\int_{|x| = R} (u+a) \nabla u \cdot n \ dS = O(R^{-1}).
\]
As a result,
\[
- \int_{\mathbb{R}^2} |\nabla u|^2dx =2\pi\gamma \int_{\mathbb{R}^2} V(u+a)^2dx.
\]
    From this we deduce that
  $\ds{\int_{\mathbb{R}^2} |\nabla u|^2dx = \int_{\mathbb{R}^2} 
(u+a)^2Vdx = 0}$.  This in
turn implies that $u \equiv 0$. But this contradicts $u=\cG(V(u+a))$. 
Hence we can not have
$Z=0$.\\

\noindent{\bf Step 4.}  It remains to show that $(u+a)Z^{-1}\ge 0$.
We only establish this when $Z>0$. The case $Z<0$ can be treated likewise.
First take a smooth function
$\varphi_{\d}: {\R}
\to (-\i,0]$ such that $\varphi'_{\d} \ge 0$ and
\[
\varphi_{\d}(r) = \begin{cases}
0 &r > -a, \\
a + r &r < -a-\d.
\end{cases}
\]
We then have
\begin{equation}
\label{eq3.11}
\int_{|x|= R}\varphi_{\d}(u)\nabla u. n dS-\int_{|x|\le R}
  \varphi'_{\d}(u)|\nabla u|^2dx = 2\pi{\gamma} \int_{|x|\le R}
\varphi_{\d}(u)(u+a)Vdx
\end{equation}
by an integration by parts. (Here $n=x/|x|$.) If $Z>0$, then we can use
\eqref{eq3.9} to assert that  $u(x)>0$ and $\varphi_\d(u(x))=0$
 whenever $|x|=R$ and R is
sufficiently large.  Since the left-hand side of
\eqref{eq3.11} is negative for such large $R$
and $(u+a)\varphi_{\d}(u) \ge 0$ we deduce
\[
\int_{\mathbb{R}^2} \varphi'_{\d}(u)|\nabla u|^2dx =
\int_{\mathbb{R}^2} V(u+a)\varphi_{\d}(u)dx = 0.
\]
We now send $\d \to 0$ to deduce
\[
0 = \int_{\mathbb{R}^2} |\nabla u|^21\!\!1(u+a \le 0)dx
  = \int_{\mathbb{R}^2} V(u+a)^21\!\!1(u+a \le 0)dx.
\]
As a result, on the set $A = \{x: a + u(x) < 0\}$ we have $\nabla u =
0$.  Hence
$u$ is constant on
each component $B$ of $A$.  But this constant can only be $-a$
because on the boundary of $A$ we
have $u+a=0$.  This is impossible unless $A$ is empty and we deduce
that $u \ge -a$ everywhere.
\qed

We now turn to the proof of Theorem 3.2. We first state and prove
a lemma. Let us write ${\Lambda}_\e$ for
$\mu^\e\big(\R^2\big)$ where $\mu^\e$ was defined right after (3.3).
\begin{lemma}\label{lem3.1}
We have that ${\Lambda}_\e>0$ for small $\e$.
Moreover
$$
\limsup_{\e\to 0}{\Lambda}_\e\le 1.
$$
\end{lemma}

{\bf Proof} Let us write $\hat u^\e$ for $(u^\e+|\log \e|){\Lambda}_\e^{-1}$ and
$\hat \mu^\e$ for ${\Lambda}_\e^{-1} \mu^\e$. By Theorem~3.1 we have that $\hat \mu^\e$
is a probability measure and $u^\e\ge 0$.
 Note that the support of the probability measure $\hat \mu^\e$ is the set
$\e K_0$. Moreover $u^\e$ is harmonic off $\e K_0$ and
\begin{equation}\label{3.17}
\hat u^\e=\cG \hat \mu^\e+\frac {|\log\e|}{{\Lambda}_\e}.
\end{equation}
 By a well-known theorem in potential theory we have
\[
\e Cap(K_0)=Cap(\e K_0)\ge \exp\left(-\frac {|\log\e|}{{\Lambda}_\e}\right),
\]
where $Cap$ denotes the logarithmic capacity.
(See for example Theorem 9.8 of \cite{Po}.)
As a result
$$
\frac {\log Cap(K_0)}{|\log\e|}-1\ge -\frac 1{{\Lambda}_\e}.
$$

 From this, we can readily deduce the  claims of the Lemma.
\qed

{\bf Proof of Theorem 3.2} It suffices to show that 
for every positive $k$,
\begin{equation}\label{eq3.15}
\lim_{\e\to 0}\sup_{|x|\le k}\left|u^{\e}(\e x)|\log \e|^{-1}+{\Lambda}_\e\right|=0,
\end{equation}
and 
\begin{equation}\label{eq3.16}
\lim_{\e\to 0}{\Lambda}_\e=\frac {\tau}{2\pi+{\tau}}.
\end{equation}

Recall that by Theorem 3.1 and Lemma 3.1
 the expression $ u^\e(\e y)|\log\e|^{-1}+1$ is nonnegative
 Also recall that  $R_0$ is 
defined so that the ball $B_{R_0}(0)$ contains the support of $V$.
Let us write $\ell_\e$ for the maximum of
$u^{\e}(\e x)$ over the ball $B_{R_0}(0)$.
We then have
\begin{eqnarray}\label{eq3.17}
u^\e(\e x)&=&\frac{\tau} {2\pi}\int \log |\e x-\e y|
\left(u^\e(\e y)|\log\e|^{-1}+1\right)V(y)dy\\
&=& {\Lambda}_\e\log \e+\frac{\tau} {2\pi}\int \log | x-y|
\left(u^\e(\e y)|\log\e|^{-1}+1\right)\ V(y)dy\nonumber\\
&\le&{\Lambda}_\e\log \e+
\frac{\tau} {2\pi}\left(\ell_\e|\log\e|^{-1}+1\right)\int_{|x-y|>1}
 \log | x-y|\ V(y)dy.\nonumber
\end{eqnarray}
Hence for every $x$ with $|x|\le k$ we have
$$
u^\e(\e x)\le {\Lambda}_\e\log \e+
c\left(\ell_\e|\log\e|^{-1}+1\right),
$$
where $c$ is a constant that  depends on $k$ only.
By choosing
$k=R_0$  we deduce that 
$$
\left (1-c|\log \e|^{-1}\right)\ell_\e\le {\Lambda}_\e \log\e +c.
$$ 
Hence 
\begin{equation}\label{eq3.18}
\ell_\e=\max_{|x|\le R_0}u^\e(\e x)\le 2{\Lambda}_\e \log\e+2c
\end{equation}
for sufficiently
small $\e$. This in turn implies
\begin{equation}\label{eq3.19}
\ell_\e |\log\e|^{-1}+1\le 2,
\end{equation}
for small $\e$. Moreover, by the second equality in (\ref{eq3.17}),
\begin{displaymath}
u^\e(\e x)|\log \e|^{-1}+{\Lambda}_\e=\frac {{\tau}}{2\pi |\log\e|}
\int\log|x-y|\left(u^\e(\e y)|\log\e|^{-1}+1\right)V(y)dy=X_1+X_2,
\end{displaymath}
where
\begin{eqnarray*}
X_1&=&\frac {{\tau}}{2\pi |\log\e|}
\int_{|x-y|\le 1}\log|x-y|\left(u^\e(\e y)|\log\e|^{-1}+1\right)V(y)dy,\\
X_2&=&\frac {{\tau}}{2\pi |\log\e|}
\int_{|x-y|\ge 1}\log|x-y|\left(u^\e(\e y)|\log\e|^{-1}+1\right)V(y)dy.
\end{eqnarray*}
Since the expression $ u^\e(\e y)|\log\e|^{-1}+1$ is nonnegative and
 bounded above for 
$y$ in the ball $B_{R_0}(0)$, we deduce that both $X_1$ and $X_2$ converge to $0$ 
in low $\e$ limit. This completes the proof of 
(3.15). 

We now turn to the proof (3.16). By the
definition of ${\Lambda}_\e$ and (3.15),
\begin{eqnarray*}
{\Lambda}_\e&=&\frac {{\tau}}{2\pi}\int 
\left(u^\e(\e y)|\log\e|^{-1}+1\right)V(y)dy\\
&=&\frac {{\tau}}{2\pi}(1-{\Lambda}_\e)\int V(y)dy+o(1)=\frac {{\tau}}{2\pi}(1-{\Lambda}_\e)+o(1).
\end{eqnarray*}
This immediately implies (3.16).
\qed

\end{section}
\bibliographystyle{plain}
\bibliography{bibcoagII}
\end{document}